\documentclass[11pt,a4paper]{amsart}
\usepackage{amssymb}
\usepackage[dvips]{graphicx}
\usepackage{psfrag}
\usepackage{amscd}
\usepackage[latin1]{inputenc}
\numberwithin{equation}{section}

\DeclareMathAlphabet{\mathpzc}{OT1}{pzc}{m}{it}
\newcommand{\bas}{{\circ\!\bigcirc\!\circ}}

\newcommand{\angl}{\measuredangle}

\def\dist{\operatorname{dist}}
\def\diam{\operatorname{diam}}

\def\ti{\tilde}

\newcommand{\cir}[1]{\overset{\circ}{#1}}
\newcommand{\gen}{\operatorname{Gen}}

\def\({\left(}
\def\){\right)}
\def\oli{\overline}

\def\raw{\rightarrow}

\def\no={\neq}
\def\sm{\setminus}

\def\mat{{\mathpzc M\mathpzc a\mathpzc t}}

\def\C{{\mathbb C}}
\def\D{{\mathbb D}}

\def\N{{\mathbb N}}

\def\Q{{\mathbb Q}}
\def\R{{\mathbb R}}

\def\T{{\mathbb T}}

\def\Z{{\mathbb Z}}

\def\AA{{\mathcal A}}
\def\BB{{\mathcal B}}
\def\CC{{\mathcal C}}
\def\DD{{\mathcal D}}

\def\FF{{\mathcal F}}

\def\II{{\mathcal I}}

\def\LL{{\mathcal L}}
\def\MM{{\mathcal M}}

\def\OO{{\mathcal O}}
\def\PP{{\mathcal P}}

\def\RR{{\mathcal R}}

\def\TT{{\mathcal T}}

\def\al{\alpha}
\def\be{\beta}
\def\ga{\gamma}
\def\de{\delta}
\def\varep{\varepsilon}

\def\th{\theta}

\def\la{\lambda}

\def\si{\sigma}

\def\De{\Delta}

\def\La{\Lambda}

\newcommand{\comm}[1]{}

\newcommand{\rt}{\operatorname{root}}

\renewcommand{\mod}{\operatorname{mod}}

\theoremstyle{plain}
\newtheorem*{Main}{Main Theorem}

\newtheorem*{Lalemma}{The $\la$-Lemma}

\newcommand{\tl}{\tilde}

\newtheorem{Thm}{Theorem}[section]
\newtheorem{Prop}[Thm]{Proposition}
\newtheorem{Lem}[Thm]{Lemma}

\newtheorem{Cor}[Thm]{Corollary}
\newtheorem{thm}{Theorem}[section]
\newtheorem{prop}[Thm]{Proposition}
\newtheorem{lem}[Thm]{Lemma}

\newtheorem{cor}[Thm]{Corollary}

\theoremstyle{remark}

\theoremstyle{definition}
\newtheorem{Def}{Definition}[section]
\newtheorem{defn}{Definition}[section]

\newcommand{\lemref}[1]{Lemma \ref{#1}}
\newcommand{\thmref}[1]{Theorem \ref{#1}}
\newcommand{\propref}[1]{Proposition \ref{#1}}

\newcommand{\circarc}[2]{[{#1}\circlearrowleft{#2}]}

\begin{document}

\title[Mating non-renormalizable quadratic polynomials]
{Mating non-renormalizable quadratic polynomials}
\author{Magnus Aspenberg and Michael Yampolsky}
\thanks{The first author was partially supported by the Foundation Blanceflor
Boncompagni-Ludovisi, n{\'e}e Bildt and by the Fields Institute.}
\thanks{The second
author was partially supported by an NSERC Discovery Grant.}

\begin{abstract}
In this paper we prove the existence and uniqueness of matings of the
basilica with any quadratic polynomial which lies outside of  the $1/2$-limb of $\MM$, is non-renormalizable, and
does not have any non-repelling periodic orbits.
\end{abstract}

\maketitle

\section{Introduction}
\label{sec:intro}
\subsection{Two definitions of mating}
The idea of mating quadratic polynomials was introduced by
Douady and Hubbard \cite{Do1} as a way to dynamically parameterize
parts of the parameter space of quadratic rational maps by pairs of
quadratic polynomials. We will present several different ways of describing the
construction, which lead to equivalent definitions in the case
which is of interest to us.

Consider two quadratic polynomials
 $f_1(z)=z^2+c_1$ and $f_2(z)=z^2+c_2$ whose Julia sets $J_1$ and
$J_2$ are connected and locally connected. For $i=1,2$ denote $\Phi_i$
the B\"ottcher coordinate at infinity
$$\Phi_i: \hat{\C} \sm K_i \raw \hat{\C} \sm \bar \D,$$ where $K_i$
is the filled Julia set of $f_i$.
It gives a conjugation
\[
\Phi_i \circ f_i (z) = (\Phi_i (z))^2, \text{  for $i=1,2$.}
\]
Carath{\'e}odory's Theorem implies that $\Phi_i^{-1}$ extends to a
continuous
parameterization $\partial\D\to J_i$. Setting
$$\ga_i: t \raw \Phi_i^{-1}(e^{2\pi it}) \in J_i,$$
we have
\begin{equation}
\label{eq:semiconj}
f_i(\ga_i(t))=\ga_i(2t).
\end{equation}

The topological space
\[
X = (K_1 \sqcup K_2) / (\ga_1(t) \sim \ga_2(-t))
\]
is obtained by glueing the two filled Julia sets along their boundaries in
reverse order.
Note that by (\ref{eq:semiconj}) the dynamics of $f_1|_{K_1}$ and $f_2|_{K_2}$
correctly defines a dynamical system $F:X\to X$,
\[
F= (f_1 \big\lvert_{K_1} \sqcup f_2
\big\lvert_{K_2}) / (\ga_1(t) \sim \ga_2(-t)).
\]

If $X$ is homeomoprhic to $S^2$, then we say that $f_1$
and $f_2$ are {\em topologically mateable}. In this case,
we call the mapping $F$ the {\it topological mating}, and
use the notation
$F=f_1 \sqcup_{\TT} f_2$.

Assume further, that there exists a homeomorphic change of coordinate
$\psi:X\to\hat\C$ which is conformal on $\cir{K_1}\cup\cir{K_2}$
and such that
$$R=\psi\circ F\circ\psi^{-1}:\hat\C\to\hat\C$$
is a rational mapping. We then say that $R$ is a {\it conformal mating}
(or simply a {\it mating}) of $f_1$ and $f_2$, and write
$R=f_1\sqcup f_2$.
The  pair of quadratics $f_1$ and $f_2$ is then called {\it
  conformally mateable}.
Conformal mateability thus implies, in particular, topological mateability.

Let us give another useful definition of mating. Let $\copyright$ be
the complex plane compactified by adjoining the circle of directions at infinity
$\{ \infty\cdot e^{2\pi i \th} : \th \in S^1\}$.
Given two quadratic polynomials $f_1$ and $f_2$ as before,
consider the extension of $f_i$ to the circle at infinity given by
$$f_i(\infty\cdot e^{2\pi i\th})=\infty\cdot e^{4\pi i\th}.$$
Glueing the two circles at infinity in reverse order, we obtain a
2-sphere $\Omega=\copyright_1 \cup \copyright_2 / \sim_{\infty}$, with the
equivalence relation $\sim_{\infty}$ identifying $(\infty\cdot e^{2 \pi i
  \th_1})$ with $(\infty\cdot e^{2 \pi i \th_2})$ whenever
$\th_1=-\th_2$, and a well defined map $f_1 \sqcup_{\FF} f_2$ equal to $f_i$
on $\copyright_i$, $i=1,2$.
The map $f_1 \sqcup_{\FF} f_2$ is called the {\em formal mating} between
$f_1$ and $f_2$.

For each $\th\in S^1$ we denote $R_i(\th)$ the
{\it external ray} of $f_i$ with angle $\th$ given by
$$\Phi_i^{-1}(\{re^{2 \pi i \th}\text{ for }r\geq 1\}).$$
Label $\hat R_i(t)$ the closure of $R_i(t)$ in $\Omega$.
We define the {\it ray equivalence relation} $\sim_r$ on $\Omega$ in the following
way: $x\sim_r y$ if and only if there exists a finite sequence of closed
external rays $\{\hat R_{i_j}(t_j)\}_{j=1,\ldots,k}$ with the property
$$\hat R_{i_j}(t_j)\cap \hat R_{i_{j+1}}(t_{j+1})\neq \emptyset,\text{
  for }1\leq j\leq k-1 \text{ and }\hat R_{i_1}(t_1)\ni x,\; \hat R_{i_k}(t_k)\ni y.$$
If $f_1$ and $f_2$ are topologically mateable then
it follows from the definition that the topological space
$\copyright_1 \sqcup \copyright_2 / \sim_{\infty}$ modulo $\sim_r$ is
again a 2-sphere and
\[
f_1 \sqcup_{\TT} f_2 = f_1 \sqcup_{\FF} f_2 / \sim_r.
\]

We can now give another equivalent definition of conformal mating in
terms of  ray equivalence: $f_1$ and $f_2$ are {\it conformally mateable}
if there exists
a rational mapping $R:\hat\C\to\hat \C$ and a pair of semiconjugacies
$\phi_i:K_i\to\hat \C,\; i=1,2$
\[
R \circ \phi_i = \phi_i \circ f_i,
\]
such that the following holds: $\phi_i$ is conformal on $\cir{K_i}$,
and $\phi_i(z)=\phi_j(w)$
if and only if $z \sim_r w$. The map $R$ is  called a
{\it conformal mating} between $f_1$ and $f_2$.

Recall that two branched coverings $F_i:S^2\to S^2$, $i=1,2$ with finite
postcritical sets $P_i$ are
equivalent in the sense of Thurston if there exist orientation preserving
homeomorphisms of the sphere $\phi$ and $\psi$ such that
$\phi\circ F_1=F_2 \circ \psi$, and $\psi$ is isotopic to $\phi$ rel $P_1$.
Using Thurston's characterization of postcritically finite rational mappings
as branched coverings (see \cite{DH}),
Tan Lei \cite{Tan} and Rees \cite{Rees1} demonstrated that if $f_i(z)=z^2+c_i$, $i=1,2$
is a pair of postcritically finite quadratics and the parameters $c_1$ and $c_2$
are not in conjugate limbs of the Mandelbrot set, then the formal mating
$f_1\sqcup_\FF f_2$ (or a certain degenerate form of it) is equivalent to a
quadratic rational map $R$ in the sense of Thurston.

Further, Rees \cite{Rees2} and Shishikura \cite{Shishikura} showed that under the
above assumptions, $f_1$ and $f_2$ are conformally mateable.

Note that the condition that $c_1$ and $c_2$ are not in conjugate limbs is
clearly necessary for topological mateability. Indeed, otherwise the
cycles of external rays $\{R_1(t_j)\}$ and $\{R_2(s_j)\}$ landing at the
dividing fixed points of the respective maps have opposite angles $t_j=-s_j$
(see e.g. \cite{Milnor2}). Thus $\{\hat R_1(t_j)\}\cup\{\hat R_2(s_j)\}$
separates $\Omega$ and therefore $\Omega/\sim_r$ is not homeomorphic to $S^2$.
It is remarkable that this condition is also sufficient when $f_1$ and $f_2$
have finite critical orbits, as this includes cases when both Julia sets
are dendrites with empty interior.

First examples of matings not based on Thurston's characterization of
rational maps appeared in the paper of Zakeri and the second author \cite{YZ}.
Before formulating it, recall that an irrational number $\theta\in(0,1)$ is of {\it bounded
type} if there exists $B>0$ such that $\theta$ can be expressed as an infinite
continued fraction with terms bounded by $B$.

\medskip
\noindent
{\bf Theorem.} {\it Let $\theta_1$ and $\theta_2$ be two irrationals of
bounded type, such that  $\theta_1 + \theta_2 \neq 1$. Then the pair of quadratic polynomials
$f_i=e^{2\pi i\theta_j}z+z^2$, $j=1,2$ are conformally mateable.}

\medskip
\noindent
The mating $R=f_1\sqcup f_2$ is unique up to a M{\"o}bius change of coordinates,
and is identified algebraically. However, it is very far from being postcritically
finite. The postcritical sets of its two critical points are quasicircles,
bounding a pair of Siegel disks. The approach taken in \cite{YZ} consists
in defining a dynamical {\it puzzle} partition of the Riemann sphere $\hat \C$ for the
mapping $R$. The renormalization theory of critical circle maps \cite{Yampolsky}
can be used to show that nested sequences of puzzle pieces shrink to points.
This provides a combinatorial description of the Julia set of $R$, sufficient
to verify that it is a mating.

The history of the problem we consider in this paper is as follows.
In 1995 J. Luo \cite{Luo} has proposed an approach to constructing a particular
class of non postcritically finite matings of the following sort.
A quadratic polynomial $f_c(z)=z^2+c$ is called {\it starlike} if
$c$ is contained in one of the hyperbolic components attached to the
main cardioid of the Mandelbrot set $\MM$. The name is due to the fact that
Hubbard trees associated to such components have only one branching point.

A {\it Yoccoz' quadratic polynomial} has only repelling periodic cycles,
and is renormalizable at most finitely many times. Yoccoz (see e.g. \cite{Hubbard})
has proved that such polynomials are combinatorially rigid, and have locally connected
Julia sets. Luo has proposed mating starlike maps with Yoccoz' ones,  arguing that
the Yoccoz' puzzle partition for quadratics can be transplanted into the
quadratic rational map. In this paper we carry this program out for a particular
instance of critically finite starlike polynomial $f_{-1}(z)=z^2-1$,
whose Julia set is known as the {\it basilica}. We use the symbol $\circ\hspace{-4pt}\bigcirc\hspace{-4pt}\circ$ as a graphical
reference to this particular quadratic parameter, to avoid awkward notation.
Thus $f_{-1}$ becomes $f_\bas$, and its Julia set is denoted $J_\bas$.
We prove:

\begin{Main} \label{main}
Suppose $c$ is a non-renormalizable parameter value outside the $1/2$-limb of $\MM$
such that $f_c$ does not have a non-repelling periodic orbit.
Then the quadratic polynomials $f_c$ and $f_\bas$ are conformally mateable,
and their mating is unique up to a M{\"o}bius coordinate change.
\end{Main}

It will be evident from the argument how to adapt it to work for an arbitrary
starlike map, however, we decided to specialize to the case $f_\bas$ for the
sake of clarity. Potentially, the methods of the proof should also work
for the case of a general Yoccoz' parameter $c$, or even an infinitely renormalizable
parameter with good combinatorics.

Since $f_\bas$ has a superattracting orbit $0 \raw -1 \raw 0$, any
candidate mating $R$ must exhibit a superattracting orbit of order
$2$. Let us place the critical point at $\infty$ and assume that
$R(\infty)=0$, $R^2(\infty)=\infty$. The following family will serve
as our candidate matings:
\[
R_a(z) = \frac{a}{z^2+2z}.
\]
The critical points of $R_a$ are $\infty$ and $-1$.

A crucial obstacle now (and a principal difference with \cite{YZ}) is that
there is no algebraic approach to specifying the candidate mating of
$f_c$ and $f_\bas$.
Instead, and similarly to Yoccoz' rigidity result, we will define a puzzle
partition in the parameter space of $R_a$, and select the mating as
the unique intersection point of a specific sequence of puzzle-pieces.

\subsection*{Acknowledgements}
We wish to thank Mitsu Shishikura for useful comments. We are grateful to Vladlen
Timorin for dicussing his own results on the family $R_a$ with us, and for a helpful
remark on a preliminary version of the paper. We are grateful to Roland Roeder for
help with computer graphics and stimulating conversations. We thank
Tan Lei for useful comments on an earlier version. Finally, the first author
thanks Rodrigo Perez for stimulating discussions.

This work was carried out at the Fields Institute, during the Dynamics Thematic Program
of 2005-6. We gratefully acknowledge the Institute's support and hospitality.

\section{Basic properties for $R_a$ and ${f_\bas}$} \label{pset}
For ease of reference, we summarize in this section some of the basic
properties
of the mapping $f_\bas(z)=z^2-1$ and the quadratic rational maps in
the family $R_a$.
We refer the reader to \cite{Milnor1} for the discussion of the
properties of Fatou
and Julia sets, and to \cite{Milnor2} for the properties of external
rays of polynomial maps.

\medskip
\noindent
\subsection{Basic properties of ${f_\bas}$}

\begin{figure} 
\begin{center}
\includegraphics[scale=.4]{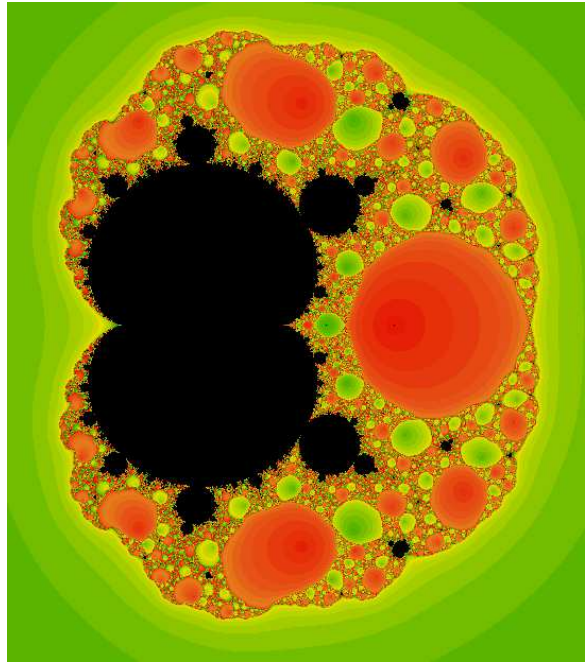}
\end{center}
\caption{The parameter set for $R_a$.\label{param} }
\end{figure}


Let us begin with the following general statement
(cf. \cite{Milnor1}).

\begin{Lem} \label{rat}
Let $U$ be a simply-connected immediate basin of  a superattracting
periodic point of a rational mapping $F:\hat\C\to\hat \C$ of period $q$.
Denote $\phi:U\mapsto \D$ a B{\"o}ttcher coordinate:
$\phi(F^q(z))=(\phi(z))^d$ for some $d>1$.
An {\it internal ray} is a curve $\phi^{-1}(\{re^{2\pi it}|\;r\in[0,1)\}.$
Then:
\begin{itemize}
\item suppose, $p$ is a repelling or parabolic periodic point on the
  boundary of $U$. Then $p$ is
the landing point of an internal ray whose period  is divisible by the
  period of $p$;
\item conversely, every periodic internal ray lands at a repelling or
  parabolic periodic
point in $\partial U$.

\end{itemize}
\end{Lem}

Let $B_0$, $B_{-1}$ be the immediate basins of attraction of $0$ and
$-1$ respectively
for $f_\bas$. Let $B_{\infty}$ be the basin of attration at
infinity. Note that
$f_\bas: B_0 \mapsto B_{-1}$ is also a $2\to 1$
covering branched at $0$.


\begin{Lem} \label{bubble}
For any two Fatou components $A$ and $B$ of $f_\bas$,  neither of which is the
attracting basin of infinity, exactly one of the following
holds:

\begin{enumerate}

\item $\oli{A} \cap \oli{B} = \emptyset$.

\item $\oli{A} \cap \oli{B}$ is only one point, which is a pre-fixed point for $f_\bas$.

\item $A=B$.

\end{enumerate}
\end{Lem}


\comm{

\begin{proof}
Clearly the Fatou components $B_0$ containing $0$ and $B_{-1}$ containing
$-1$ touch at only one point, since $f_\bas^2$ is a double covering
from $\oli{B_0}$ to itself, we have that $\partial B_0$ contains a
fixed point of $f_\bas^2$. But since $0,-1$ are fixed points
for $f_\bas^2$ and not for $f_\bas$, there must be a fixed point $x$ for
$f_\bas$ on the boundary of $B_0$. Since $x$ is also a fixed point on the
boundary of $B_{-1}$, we have that $B_0$ and $B_{-1}$ touch at $x$.

Since $f_c^2$ maps $B_0$ (and $B_{-1}$) onto itself, we can apply
Lemma \ref{rat} to $f_c^2$ and conclude that there is a internal ray
starting at $0$ and landing at the fixed point $x$. Its image is an
internal ray in $B_{-1}$ landing at $x$ as well.

There cannot be another fixed point
for $f_\bas$ on the boundary of $B_0$, since the internal landing ray
landing at $x$ from Lemma \ref{rat} has
to have internal angle $0$, since it is fixed and there is no other
internal ray which is fixed under the map $d(\th)=2\th$.

Since every Fatou component is a preimage $B_0$ the claim follows.
\end{proof}

}

\noindent
The statement of the Lemma follows immediately from the Maximum Principle.
Note, that the boundaries of the Fatou components $B_0$ and $B_{-1}$
touch at the repelling fixed point $\alpha$ of $f_\bas$.


Since the mapping $f_\bas$ is hyperbolic, its Julia set is locally
connected. In particular, if $\Phi: \hat{\C} \sm K(f_\bas)\mapsto \C\sm\D$
denotes the B{\"o}ttcher coordinate at $\infty$, the Carath{\'e}odory's
Theorem implies that
$\Phi^{-1}$ extends continuously to $\partial\D$. Moreover, every external
ray $R(\theta)=\Phi^{-1}(\{re^{2\pi i\theta}|\;r>1\})$ lands at
a point of the Julia set. We denote
\[
\ga(\th) = \lim_{r \raw 1^+} \Phi(r e^{2 \pi i \th} ).
\]


Hyperbolicity of $f_\bas$ also implies:

\begin{Lem} \label{shrink}
Let $F_i$ be an arbitrary infinite sequence of distinct Fatou components
of $f_\bas$. Then $\diam F_n\to 0$.
\end{Lem}

\noindent
We will also make use of the following Lemma:

\begin{Lem} \label{pminus1}
A point $z \in J_\bas$ is a landing point
of precisely two external rays if and only if $z$ is a preimage of the
fixed point $\al$. No other point $z\in J_\bas$ is biaccessible.
\end{Lem}

\noindent
The angles of the two external rays which land at $\alpha$ are easily
identified as $1/3$ and $2/3$.

\comm{
\begin{proof}
The $\al$ fixed point $z_0 = 1/2 - \sqrt{5}$ has two landing
rays with angles $1/3$ and $2/3$. Indeed, since $z_0$ is the common boundary of
$B_0$ and $B_{-1}$, $z_0$ divides the filled in Julia set
$K_{-1}=K(f_\bas)$ into two connected parts, i.e. $K_{-1} \sm z_0$ is
two connected components. Therefore $z_0$ is the landing point of
precisely two external rays, namely those with angles $1/3$ and $2/3$.

Since every point $z$ on the Julia set for
$f_\bas$ is the boundary of the union of all preimages of $B_0$, we
have that every preimage to $z_0$ has exactly two landing rays.

By construction, these preimages are exactly those points which have
the property that if we delete one of these points from $K_{-1}$, the
remaing set will be two connected subsets of $K_{-1}$.
\end{proof}
}

\subsection{Properties of maps in the family $R_a$}
In what follows, we will refer to the illustration of
the parameter space for the family $R_a$  pictured in Figure \ref{param}.

For $R_a$ let $A_{\infty}$ be the immediate basin of attraction at infinity,
and $A_0$ the Fatou component containing $0$.

\begin{figure}[h] 
\begin{center}
\includegraphics[scale=.4]{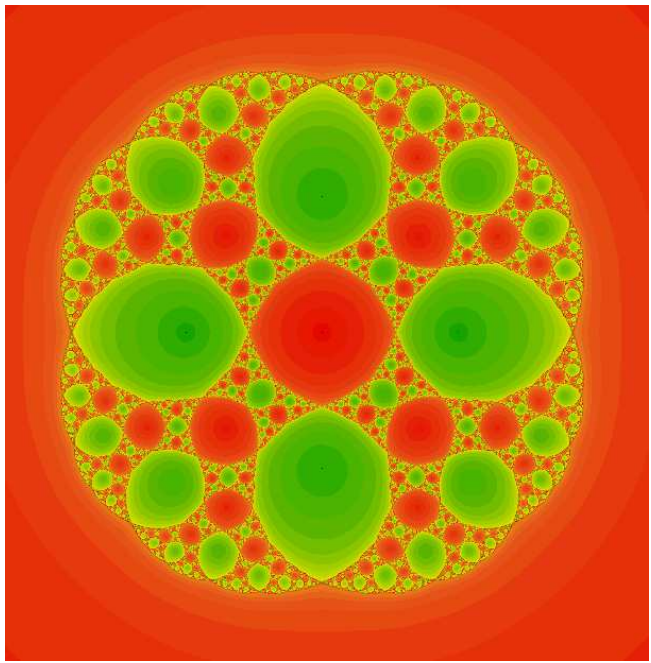}
\end{center}
\caption{A capture dynamics: dynamical plane of $R_{2}$. \label{cross}}
\end{figure}

\begin{figure}[h] 
\begin{center}
\includegraphics[scale=.4]{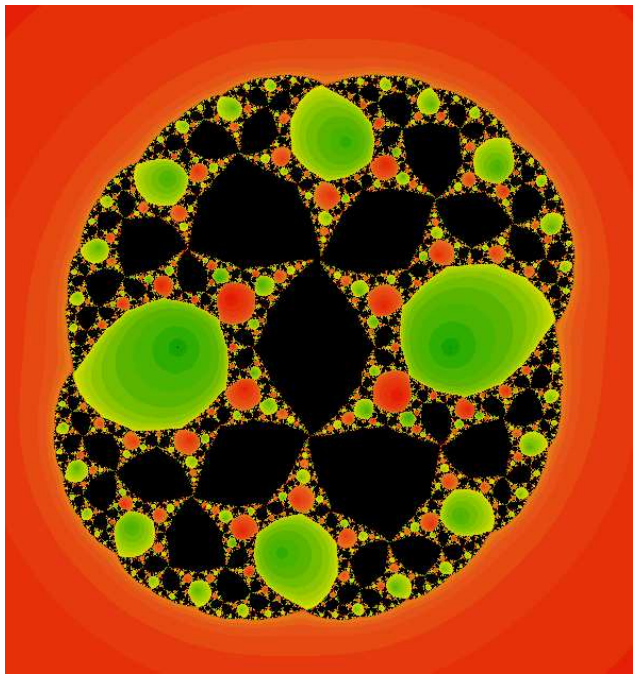}
\end{center}
\caption{Stony Brook preprint cover: the dynamical plane of the mating of basilica and  Douady's rabbit. \label{rabbit}}
\end{figure}

Let us note:
\begin{prop}
\label{first bubbles}
The Fatou components $A_0$ and $A_\infty$ are distinct and simply-connected.
The critical point $-1$ of $R_a$ is never contained in $A_\infty$.
\end{prop}
\begin{proof}
We have $A_0\neq A_\infty$ by Denjoy-Wolff Theorem.
%
If $A_\infty$ is multiply-connected, then, necessarily, $-1\in A_\infty$, by the Riemann-Hurwitz formula.
Thus $A_0$ contains all critical values of $R_a$. In this case, it follows
(see e.g. \cite{Milnor-remarks}, Lemma 8.1) that the Julia set of $R_a$ is totally disconnected, and that every orbit in the Fatou set converges to an attracting fixed point, which is impossible.
\end{proof}

Note, that whenever $a$ is such that $-1\in A_0$, the Fatou set of $R_a$
is the union of $A_0$ and $A_\infty$. The Julia set $J(R_a)$ is the common
boundary of the two Fatou components, and we have (see, for instance,
\cite{CG}, Theorem 2.1 on p. 102):

\begin{prop}
If $-1\in A_0$, then $J(R_a)$ is a quasicircle.
\end{prop}

\noindent
In the parameter space (Figure \ref{param}) the above values of $a$ form the
``exterior'' hyperbolic component which we denote $P_\infty$.

More generally, a {\it capture} hyperbolic component for the family
$R_a$ contains maps for which there exists
an iterate $R_a^{n}(-1)\in A_\infty$.
The smallest such $n$ will be referred to as the {\it generation} of the capture
component.

For instance, $a=2$ is the center of the biggest red ``bubble'' in Figure
\ref{param},
in which we have $R_a^2(-1)\in A_\infty$. The corresponding Julia set
is depicted in Figure \ref{cross}.

\comm{

In this case there is a conjugacy
between $\FF_{hyp}$ and $\FF(p) \sm B_{\infty}$, where $\FF_{hyp}$ is
the union of all preimages of $A_{\infty}$. Indeed,
there is a Riemann mapping $h$ mapping
$A_{\infty}$ onto $B_0$. For any $z \in A_{\infty}$, for $w=R_a(z)$,
we can extend $h$ to $A_0$ by setting
\begin{equation} \label{conj}
h(w) = h \circ R_a (z) = f_\bas \circ h(z).
\end{equation}
Since there is no critical point in any preimage of neither
$A_{\infty}$ under $R_a$ nor $B_0$ under $f_\bas$ we can extend $h$
to the entire set $\FF_{hyp}$.

If the critical value $-a$ belongs to some preimage $P$ of $A_{\infty}$
then the Fatou set $\FF_a$ is homeomorphic to $\FF(p) \sm B_{\infty} /
\sim$, where $\sim$ is an equivalence relation defined in the following way:
Again let $h$ be the Riemann mapping above, mapping $A_{\infty}$ onto $B_0$.
We can extend $h$ by the functional equation (\ref{conj}) until we hit
a component $P$ which contains the critical value $-a$ fr $R_a$.
At this point the preimage of $P$ under $R_a$ is only one component
whereas the preimage of $h(P)$ under $f_\bas$ consists of two
components $P_1$ and $P_2$. By identifying these two components we can
still extend $h$ to $\FF_a$, but its image will be a
subset of $\FF(p) \sm B_{\infty}$. To put it in another way, we get a
homeomorphism $h:\FF_a \mapsto \FF(p)\sm B_{\infty} / \sim$, where
$\sim$ is an equivalence relation defined by identifying all preimages
of $P_1$ and $P_2$.

}

\noindent
Similarly to the statement of Lemma \ref{bubble}, we will show in \S\ref{sec:parabubble}:

\begin{Lem} \label{Ra}
Suppose that the parameter $a$ is chosen outside of the closure $\bar P_\infty$.
Then
given any two Fatou components $A$ and $B$ in the basin of $\infty$ of  $R_a$ exactly one of
the following holds:

\begin{enumerate}

\item $\oli{A} \cap \oli{B} = \emptyset$,

\item $\oli{A} \cap \oli{B}$ is only one point,

\item $A=B$.

\end{enumerate}
Moreover, if the case (2) occurs, then $\bar A\cap \bar B$ is either a preimage of the
fixed point $$x_a\equiv \bar A_0\cap\bar A_\infty$$ or a pre-critical point.
For the latter possibility to occur, the parameter $a$ must belong to the boundary
of a capture component.

\end{Lem}

Denote $\mat$ the set of parameter values $a$ not contained in any of
the capture
components. This set is colored in black in Figure \ref{param}.
The interior of $\mat$ contains matings with basilica, and thus should
be naturally identified with $\cir{\MM}$ with the $1/2$-limb removed.

As an example of a mating in $\mat$, consider Figure \ref{rabbit}.
This image was popularized on the cover of Stony Brook preprint series;
it is the mating of Douady's rabbit with basilica.

\section{Orbit portraits for quadratic polynomials} \label{orbit}
In this section we provide a brief summary of several results on the
combinatorics of external rays of quadratic polynomials  following Milnor's paper
\cite{Milnor2}. All proofs are given in \cite{Milnor2}.

Let the points $\{x_1,x_2=f(x_1), \ldots,x_p=f(x_{p-1})\}$ form a periodic orbit of a
quadratic polynomial $f_c(z)=z^2+c$ with period $p$.
Assume further, that this orbit is either repelling or parabolic,
and hence the landing set of a finite collection of periodic external rays $R(\theta_i)$
(see e.g. \cite{Milnor1}).

\begin{defn}
For each $1\leq i\leq p$ let $A_i=\{\theta^i_1,\ldots,\theta^i_k\}$ denote the set of angles of the
external rays landing at $x_i$. The collection $\OO=\{A_1,\ldots, A_p\}$ is called the {\em orbit}
portrait of the cycle $(x_1,\ldots,x_p)$. According to the type of the cycle, the orbit portrait
is either {\it repelling} or {\it parabolic}.
\end{defn}

\noindent
Given the periodicity of $x_i$, the iterate $f_c^i$ permutes the rays with angles in $A_i$.
The following is immediate:

\begin{Lem} \label{basicorb}
Given an orbit portrait $\OO=\{A_1,\ldots,A_p\}$ the size of $A_i$ is the same for all $i$.
Moreover, $A_{i+1}=2 A_i \mod\Z$, and if $|A_i|\geq 3$, then
the cyclic order of the angles $\theta^i_j\in A_i$ is the same as that of
their images $2\theta^i_j\mod\Z\in A_{i+1}$.
\end{Lem}

\begin{defn}
For $A = \{\th_1,\ldots,\th_k\}\subset \T$, write 
$\exp(A)= \{ e^{ 2 \pi i \th_1}, \ldots , e^{2 \pi i \th_k} \}\subset S^1$. A {\it formal orbit portrait} is a collection $\{A_1,\ldots,A_p\}$ of subsets of $\T$ for which
the following properties hold:
\begin{itemize}
\item each $A_i$ is a finite subset of $\T$;
\item for each $j$ modulo $p$, the doubling map $t\mapsto 2t\mod \Z$ carries $A_j$ bijectively onto $A_{j+1}$
preserving the cyclic order around the circle;
\item all of the angles in $A_1\cup\cdots\cup A_p$ are periodic under doubling with the same period $rp$;
\item for each $i\neq j$, the convex hulls of the sets $\exp(A_i)$ and $\exp(A_j)$ are disjoint.
\end{itemize}
\end{defn}

\noindent
The {\em valence} of an orbit protrait $\OO$ is $v_\OO=|A_i|$.
Every angle in $A_i$ is periodic of period $pr$. Since there
are $pv_\OO$ angles in $\OO$, the quantity $v_\OO/r$ is the
number of distinct cycles of external rays in the orbit portrait $\OO$.

\begin{lem}
Only two possibilities can occur: either $v_\OO=r$ or $v_\OO=2$
and $r=1$.
\end{lem}

\noindent
Assume that $v_\OO\geq 2$.
For each $A_i$, the complement $\T\setminus A_i$ consists of finitely many
{\em complementary arcs}. Each such arc corresponds to a sector between two
of the rays landing at $x_i$.

\begin{Lem} \label{critval}
Let $\OO=\{A_1,\ldots,A_p\}$ be a formal orbit portrait. Then
every complementary arc for $A_i$, except for one is mapped
one-to-one under $z\mapsto 2z$ onto a complementary arc of $A_{i+1}$.
The exception is the  {\em critical arc} of $A_i$, which has length
greater than $1/2$.  The  image of the critical arc wraps around the whole unit
circle, covering one of the complementary arcs of $A_{i+1}$ twice.

If the portrait $\OO$ is realized by a quadratic polynomial, then
for each $i$, the sector corresponding to the critical arc of $A_i$
contains the critical point $0$.
\end{Lem}


\begin{Lem}
Assume that $v_\OO\geq 2$. There exists a unique shortest complementary
arc in $\OO$.
If the portrait is realized by a quadratic polynomial $f_c$, then
the sector corresponding to this arc can be characterized
among the $pv$ sectors formed by the rays landing at points $x_i$
as the one which contains the critical value $c=f_c(0)$ and no
points of the orbit $x_i$.
\end{Lem}

\begin{defn}
The complementary arc in the previous lemma is referred to as the {\it characteristic arc} of the orbit
portrait.
\end{defn}



\section{Bubble rays} \label{raysandaxes}

To construct a Yoccoz puzzle partition for the quadratic rational maps in $\mat$, we
will use chains of Fatou components in place of external rays. This method was
employed in \cite{YZ} and \cite{PR1}, it was also suggested in \cite{Luo}. We begin by
describing such chains in the filled Julia set of $f_\bas$; this discussion,
while mostly trivial, will serve as a useful preparation for handling
maps in the family $R_a$.


\subsection{Bubble rays for $f_\bas$}

Recall that $B_0$ and $B_{-1}$ denote the components of the immediate
super-attracting basin of $f_\bas$, labelled according to the point
in the critical orbit they surround.

\begin{defn}
\label{bubble-basilica}
A {\it bubble} of $K_\bas$ is a Fatou component $F\subset \cir{K}_\bas$.
The {\it generation} of a bubble $F$ is the smallest non-negative  $n=\gen(F)$
for which $f_\bas^n(F)=B_0$.
The {\it center} of a bubble $F$ is the preimage $f_\bas^{-\gen(F)}(0)\cap F$.

If $F\neq B_0$, then let $G$ be the bubble with the lowest value of $\gen(G)$ for
which $\bar G\cap \bar F\neq \emptyset.$ We will refer to $G$ as the {\it predecessor}
of $F$, and to the point $x=\rt(F)\equiv \bar G\cap\bar F$ as the {\it root} of $F$.

A {bubble ray} $\BB$ is a collection
of bubbles $\cup_{0}^{m\leq \infty}F_k$ such that for each $k$
the intersection $\oli{F_k}\cap\oli{F_{k+1}}=\{x_k\}$ is a single point, and
$\gen(F_k)<\gen(F_{k+1})$.
\end{defn}

\noindent
Note that by \lemref{bubble}, each of the points $x_k$ is a preimage of the
$\alpha$-fixed point of $f_\bas$.
If $m<\infty$, we will refer to the component $F_m$ as the {\it last bubble} of $\BB$.
Hyperbolicity of $f_\bas$ readily implies:

\begin{prop}
There exist $s\in(0,1)$, and $C>0$ such that for a bubble $F\subset \cir{K}_\bas$
we have $$\diam (F)\leq Cs^{\gen(F)}.$$
In particular, for each infinite bubble ray $\BB=\cup_0^\infty F_k$ there exists
a unique point $x\in J_\bas$ such that $F_k\to x$ in Hausdorff sense.
\end{prop}

\noindent
We refer to $x$ as the {\it landing point} of $\BB$. By \lemref{bubble}
we have:

\begin{prop}
\label{single landing point}
If two bubble rays $\BB_1$, $\BB_2$ have the same landing point, then
one of them is contained in the other one.
\end{prop}

\noindent
By \lemref{rat}, each pre-periodic point on the boundary of a bubble is a landing point of an internal
ray. We may therefore define:

\begin{defn}
The {\it axis} of a bubble ray $\BB=\{F_k\}_0^{m\leq \infty}$ is the closed union $$\gamma(\BB)\equiv\oli{\cup_0^m\gamma_k},$$
where $\gamma_k$ for $k\geq 1$ is the union of two internal rays of $F_k$ connecting its center to the
points $x_{k-1}$ and $x_k$, and $\gamma_0$ is the internal ray of $F_0$ terminating at $x_0$.
\end{defn}

\noindent
Let $x$ be the landing point of an infinite bubble ray $\BB$. As the Julia set $J_\bas$ is
locally connected, Carath{\'e}odory's Theorem implies that there exists at least one external
ray $R(\theta)$ landing there. By \lemref{pminus1}, such $\theta$ is unique. Let us refer to
the number $-\theta$ as the {\it angle} of the bubble ray $\BB$ and denote it
$$\angl(\BB)\equiv -\theta.$$
By \propref{single landing point}, $\angl(\BB_1)=\angl(\BB_2)$ implies that one of these
rays is a subset of the other.

\noindent
We will call a bubble ray $\BB$ {\it periodic} if the angle $\angl(\BB)$ is periodic under
 doubling; the {\it period} of the ray will refer to the period of its angle.

\noindent
Note that the angle of a bubble ray can be determined intrinsically, from the
choice of the bubbles themselves. Indeed, consider the {\it spine}
$$\ell(K_\bas)\equiv K_\bas\cap \R=[-\beta,\beta],$$
where $\beta$ is the non-dividing fixed point of $f_\bas$.
The spine may also be seen as the union of the axes of the bubble rays $\BB_+$, $\BB_-$
starting with the bubble $B_0$ and terminating at $\pm\beta$ respectively.

Let $\BB=\cup F_k$ be an infinite bubble ray, landing at $x\neq \beta$. Consider the forward iterates
$x_k=f_\bas^k(x)$. Define a sequence $s(\BB)=(s_i)_1^\infty$ of $0$'s and $1$'s as follows.
We set
\begin{itemize}
\item $s_i=0$ if $x_i$ is above the spine, or equivalently, if there is a bubble $F_k$ with
$k\geq i$ which is above the spine;
\item $s_i=1$ if $x_i$ is below the spine, or equivalently, if there is a bubble $F_k$ with
$k\geq i$ which is below the spine;
\item if $i$ is the first instance when neither of these two possibilities holds,
set $s_i=1$, and $s_j=0$ for all $j>i$ (note, that in this case we necessarily have $x_i=-\beta$.

\end{itemize}
\noindent
For $\BB\subset \BB_+$ we set $s(\BB)=(0)_0^\infty.$

\noindent
We will sometimes refer to the dyadic sequence $s(\BB)$ as the {\it intrinsic address} of $\BB$.
Noting that $$(\beta,+\infty)=R(0),\text{ and }(-\infty,-\beta)=R(1/2),$$
we immediately have

\begin{prop}
For each infinite bubble ray $\BB$ we have
$$\angl(\BB)=\displaystyle-\sum_{i=1}^\infty 2^{-i}s_i, \text{ where }s(\BB)=(s_i)_0^\infty.$$
\end{prop}

\begin{figure}

\centerline{\includegraphics[height=0.8\textheight]{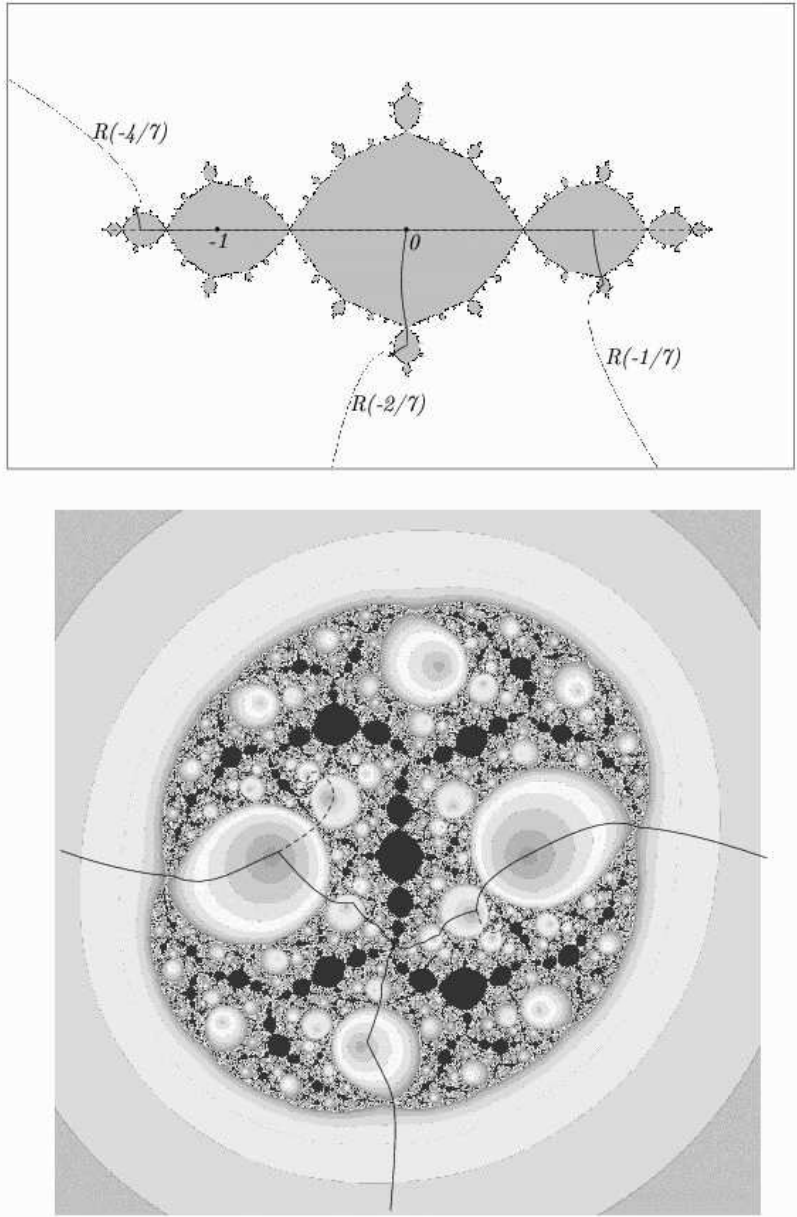}}
\caption{Bubble rays for $f_\bas$ and $R_a$. The picture below is a mating of $f_\bas$ with
a hyperbolic parameter in the $1/3$-limb of $\MM$. Three periodic bubble rays land at a
repelling fixed point of the rational map. The solid lines follow their axes. Their angles
are $1/7$, $2/7$, and $4/7$ respectively. The axes of the same bubbles are shown inside
$K_\bas$ in the above pictures. The broken lines show the position of the spines.\label{bubble rays figure}} 
\end{figure}

\subsection{Bubble rays for a map $R_a$}
The definition of a bubble ray for a rational mapping $R_a$ is completely
analogous to Definition \ref{bubble-basilica}.
\begin{defn}
\label{bubble-rat}
A {\it bubble} of $R_a$ is a Fatou component $F\subset \cup R_a^{-k}(A_\infty)$.
The {\it generation} of a bubble $F$ is the smallest non-negative  $n=\gen(F)$
for which $R_a^n(F)=A_\infty$.
The {\it center} of a bubble $F$ is the preimage $R_a^{-\gen(F)}(\infty)\cap F$.

A {bubble ray} $\BB$ is a collection
of bubbles $\cup_{0}^{m\leq \infty}F_k$ such that for each $k$
the intersection $\oli{F_k}\cap\oli{F_{k+1}}=\{x_k\}$ is a single point, and
$\gen(F_k)<\gen(F_{k+1})$.

\end{defn}
The structure of bubble rays for $R_a$ is particularly easy to describe when
$a\in\mat$, and somewhat more difficult in the capture case. We consider the
simpler possibility first.

\medskip
\noindent
{\bf The case $a\in\mat$}
Consider the B{\"o}ttcher coordinates $b_1:\DD\to B_0$, and $b_2:\DD\to A_\infty$.
The identification
$$
\phi\equiv b_2\circ b_1^{-1}:B_0\to A_\infty
$$
conjugates the dynamics of $f_\bas$ and $R_a$.
Note that by Lemmas \ref{rat} and \ref{Ra} the components $A_\infty$ and $A_0$
have a single common boundary point $x=\lim_{r\to 1-}b_2(r)$ and is fixed by
the dynamics of $R_a$. By \lemref{Ra} we have the following:

\begin{prop}
\label{touch points}
If two bubbles $F_1$ and $F_2$ of $R_a$ touch at a boundary point $z$, then
$z$ is a preimage of $x$.
\end{prop}

\noindent
By \lemref{rat}, the axis $\gamma(\BB)$ of a bubble ray $\BB$ of $R_a$ can be defined as before.
To define the {\it spine} $\ell_a$ begin by considering the union of internal rays
$l_\infty\subset\hat\C$ which is the image under $b_2$ of the segment $(-1,1)$. Let $l_0\subset A_0$ be its preimage, and set
$$t_1=\bar l_\infty\cup \bar l_0\cup\{x\}.$$
We now inductively define $t_n=t_{n-1}\cup h_1\cup h_2$ where $h_i$ are the two components of $R_a^{-1}(t_{n-1}) \sm A_\infty$ intersecting $t_{n-1}$.

\begin{defn}
We set $$\ell_a=\cup t_n,$$
and endow this arc with positive orientation as induced by the orientation of
$(-1,1)\mapsto l_\infty$.
Further, for a bubble $F$ of $R_a$ with $F\cap\ell_a=\emptyset$, we say that $F$ is
{\it above } the spine, if the unique finite bubble ray connecting it to the spine
lies above $\ell_a$ with respect to the orientation of $\ell_a$. In the complementary case,
we say that the bubble $F$ is {\it below} the spine.
\end{defn}

We define the {\it intrinsic address} $s(\BB)$ of a bubble ray $\BB$ in exactly the same fashion as
before.

The oriented spine allows us to extend inductively the conjugacy $\phi:f_\bas^{-n}(B_0)\to R_a^{-n}(A_\infty)$
so that:

\begin{prop}
\label{marking-mated}
Denote $$L=\cir{K}_\bas\cup \left(\bigcup_{n=0}^\infty f_\bas^{-n}(\alpha)\right).$$
Then $\phi$ extends as a conjugacy to the whole of $L$. Moreover,
this conjugacy obeys the property:
$$s(\phi(\BB))=s(\BB)$$
for each bubble ray $\BB$ in $K_\bas$.
\end{prop}

\begin{defn}
For an infinite bubble ray $\BB$ of $R_a$ we set the {\it angle of } $\BB$ equal to
$$\angl(B) \equiv \angl(\phi^{-1}(\BB)).$$
\end{defn}

\noindent
By construction, we have
\begin{equation}
\label{angle-address}
\angl(\BB)=\sum_{n=1}^\infty 2^{-n}s_n,\text{ where }s(\BB)=(s_n)_1^\infty
\end{equation}
for each bubble ray $\BB$ of $R_a$.

\medskip
\noindent
{\bf The case when $a$ belongs to a capture component.}
Let us exclude the trivial possibility when the critical value $-a=R_a(-1)\in A_\infty$, and denote
$n>1$ the smallest natural number for which $R_a^n(-1)\in A_\infty$ holds.
The conjugacy $\phi$ can still be extended consistently with the orientation
to $f_\bas^{-(n-1)}(B_0)$. Denote $F\ni -a$ the bubble  of $R_a$ containing the critical value, and set $H=\phi^{-1}(F)\subset \cir{K}_\bas.$

\begin{defn}
We define an equivalence relation $\sim$ on $\cir{K}_\bas$ as follows.
Connect the two preimages $H_1$, $H_2$ of $H$ by a simple arc $h\subset \hat\C\setminus K_\bas$.
 The equivalence relation identifies any two bubbles $G_1$, $G_2\subset \cir{K}_\bas$
if there exists $l\geq 0$ such that $G_1$ is connected to $G_2$ by a component of $f_\bas^{-l}(h)$.
For points $x_i\in G_i$ we set $x_1\sim x_2$ if this happens, and if
$f_\bas^{n+l}(x_1)=f_\bas^{n+l}(x_2).$
\end{defn}

\noindent
One readily verifies:
\begin{lem}
\label{capture conj}
In the capture case, the mapping $\phi$ extends as a surjective conjugacy from
$$\left.\left(\cir{K}_\bas\cup \bigcup_{i=0}^\infty f_\bas^{-i}\alpha\right)\right/_{z_1\sim z_2}\longrightarrow\bigcup_{i=0}^\infty R_a^{-i}(A_\infty\cup \{x\}).$$
\end{lem}

\section{Parabubble rays.}
\label{sec:parabubble}

Removing the $\alpha$-fixed point from the basilica $K_{\bas}$
separates it into two connected components. We will denote them
$\LL$ for ``left'', and $\RR$ for ``right''.
Put $\RR_e = \RR \sm \oli{B}_0$,
(the subscript $e$, standing for "exterior" of the right half of the
basilica). As we will see below, there is a natural correspondence
between the components of the interior of $\RR_e$, and the capture
hyperbolic components in the parameter plane of the family $R_a$.

For the remainder of this section, let us fix the notation $R_a(z)=R(z,a)$, $R_a^n(z)=R^n(z,a)$.

\begin{defn}
Let $a_0$ be such that $R_{a_0}^k(-1)=\infty$ for some $k\in\N$, and let
$n$ be the smallest such value of $k$.
Then a connected set $P$ of parameters $a$ containing $a_0$,
such that $R^n(-1,a) \in A_{\infty}$ is called a {\it capture} hyperbolic
component or a {\em parabubble}.
The point $a_0$ is called a {\it center} of $P$.
We will see further that it is unique.

Finally, we say that the
{\it generation} of $P$ is $n$, and write $\gen(P)=n$.
\end{defn}

Set $\xi_n(a)=R^n(-1,a)$. Then we have
\begin{equation} \label{xien}
\xi_{n+1}(a)=\frac{a}{(\xi_n(a))^2+2\xi_n(a)}=\frac{a}{\xi_n(a)(\xi_n(a)+2)}.
\end{equation}

\comm{
\begin{Lem} \label{xin-lemma}
Let $n \geq 2$. Then the following holds.
\begin{itemize}
\item The function $\xi_n(a)$ has simple poles at centers of capture components of
generation equal to $n$.
\item All critical points of $\xi_n(a)$ are
precisely the centers of capture components of lesser generation than
$n$ but greater than $2$,
and coincide with the poles and zeros of $\xi_n(a)$ of higher multiplicity
than $1$.
\end{itemize}
\end{Lem}

\begin{proof}
We prove the assertions of the lemma by induction. They are clearly true for $n=2$.
Assume now that they hold for some  $n \geq 2$. We have
\begin{equation} \label{xien}
\xi_{n+1}(a)=\frac{a}{(\xi_n(a))^2+2\xi_n(a)}=\frac{a}{\xi_n(a)(\xi_n(a)+2)}.
\end{equation}
We know that $\xi_n(a)=-2$ has only simple roots, since by the induction
assumption $\xi_n'(a) \neq 0$ at such a root. Hence, $\xi_{n+1}(a)$ has
simple poles at centers of all capture components of generation equal to
$n+1$. Moreover, by (\ref{xien}) all other zeros and poles of
$\xi_{n+1}(a)$ are the zeros and poles of $\xi_n(a)$ interchanged, and
for odd $n \geq 3$ we have $a=0$ is a simple zero.
Moreover, the simple poles of $\xi_n$ now
become critical points of $\xi_{n+1}$ by (\ref{xien}).
\end{proof}
}
From (\ref{xien}) it follows by a straightforward induction, that

\begin{lem}
\label{xi-count}
The degree of $\xi_n$ is the nearest integer value to $2^{n+1}/3$.
\end{lem}

We now state:
\begin{lem}
\label{number}
For $n \geq 2$, the degree of $\xi_n$ is equal to the number of
bubbles of generation $n$ in the basilica which are contained in $\RR$.
\end{lem}
\begin{proof}
To each bubble $B\subset K_\bas$ we
associate an interval
$(a,b)=I_B\subset \R/\Z$, where $a,b$ are the angles of the external rays meeting at
the root of $B$.
It is easy to see that the centers of the intervals $I_B$ of all bubbles of
generation $n$ are symmetrically distributed around the unit circle and that
each $I_B$ does not intersect $1/3$ or $-1/3$. It is easy to verify
that the
closest integer to $2^n \times (2/3)$ is equal to the number of $I_B$ which
are contained in the interval $(-1/3,1/3)$. The claim follows from \lemref{xi-count}.
\end{proof}

Denote $A_{\infty}^a$ the set $A_{\infty}$ for the map $R_a$. Let
$$\Phi_a:A_{\infty}^a \mapsto \hat{\C} \sm \D$$
 be the B\"ottcher coordinate for
$R_a$ normalized so that $\Phi_a'(\infty)>0$.
Note that $\Phi_a$ is analytic in $a$.
A direct calculation implies
\begin{equation}\label{phi}
\Phi_a(z) = \frac{1}{2}z + o (1), \text{as $z \raw \infty$}.
\end{equation}
If $-a \in A_{\infty}$ then $\Phi_a$ can be extended around $\infty$
until we hit a critical point $z=1 \pm \sqrt{1-a}$ for $R_a^2$. However,
the Green's function $g(z,a)=\log |\Phi_a(z)|$ is still well defined on
$A_{\infty}$ and moves continuously with $a$, and $g(z,a) \raw 0$ as
$z \raw \partial A_{\infty}$ for all $a \in \C$.
Let $P_{\infty}$ be the open set of parameters where $-a \in
A_{\infty}$, that is, where $J(R_a)$ is a quasicircle.
This capture component obviously contains an open neighborhood of $\infty$.
%

By the $\la$-Lemma of \cite{MSS} we have:
\begin{Lem}
The Julia set $J(R_a)$ moves holomorphically for all $a \in
P_{\infty}$.
\end{Lem}

\noindent
Let us continuously extend the Green's function $g(z,a)$ on the whole
sphere so $g(z,a)=0$ outside $A_{\infty}$. The proof of Theorems III.3.2
in \cite{CG} can be easlily adjusted to the family
$R_a^2: A_{\infty} \mapsto A_{\infty}$, to show that the Green's
function $g$ is uniformly H\"older $\al$-continuous for $|a|
\leq C$, some $\al=\al(C) \in (0,1]$. As a consequence, $g(-a,a) \raw
  0$ as $a \raw \partial P_{\infty}$, (see Theorem III.3.3 \cite{CG}).
Moreover, by (\ref{phi}), the function $\Phi_a(-a)$ has a
simple pole at $\infty$. Since $g(-a,a) \raw 0$ as $a \raw \partial
P_{\infty}$, the Argument Principle implies that $\Phi_a(-a)$ takes every
value in $\hat{\C} \sm \D$ exactly once. We get the following:
\begin{Lem}
The set $P_{\infty} \cup \{\infty \}$ is simply connected and
$P_{\infty}^c$ has logarithmic capacity equal to $1/2$.
\end{Lem}

It is easy to verify that $A_{\infty}$ does not necessarily move
continuously at $\partial P_{\infty}$ if we step inside $P_{\infty}$
(e.g. at $a=3$), but the following holds.
\begin{Lem} \label{above}
The set $\oli{A}_{\infty}$ moves holomorphically for all parameters
$a \in (\oli{P}_{\infty})^c$.
We have $a \in \partial P_{\infty}$ if $-a \in \partial
A_{\infty}^a$.
\end{Lem}
\begin{proof}
Put $\psi_a=\Phi_a \circ \Phi_{a_0}^{-1}$. Then $\psi_a$
maps $A_{\infty}^{a_0}$ onto
$A_{\infty}^a$.
If $a \notin \oli{P}_{\infty}$ then $-a
\notin A_{\infty}$ by definition and we have that
$\psi_a(z)=\psi(z,a)$ is a holomorphic motion on
$A_{\infty}^{a_0} \times \oli{P}_{\infty}^c$. By the $\La$-Lemma,
$\partial A_{\infty}^a$ also moves holomorphically.

If $-a_1 \in \partial A_{\infty}$ for
some $a_1 \notin \oli{P}_{\infty}$ then since $A_{\infty}$ moves
holomorphically, the point
$-a_1$ is an image of some point $z_1 \in \partial A_{\infty}^{a_0}$ under
$\psi$, i.e. $\psi(z_1,a_1) = -a_1$. The analytic function $\psi_{z_1}(a)$
satisfies $\psi_{z_1}(a_1)+a_1=0$. Either
$\psi_{z_1}(a) + a \equiv 0$ or not. If so, then $-a \in \partial A_{\infty}$
for all $a \in (\oli{P}_{\infty})^c$, which is clearly false. If not so,
then choose a small disk $B(a_1,\varep) \subset (\oli{P}_{\infty})^c$
and some $z_2 \in A_{\infty}^{a_0}$, with $z_2$ sufficiently
close to $z_1$, such that $|\psi_{z_1}(a)-\psi_{z_2}(a)| <
|\psi_{z_1}(a)+a|$ for $a \in \partial B(a_1,\varep)$.
By Roche's Theorem,
$\psi_{z_2}(a)+a=0$ must have a solution $b \in B(a_1,\varep)$, which
means that $-b \in A_{\infty}^b$, which is a contradiction.

\comm{

Moreover, $A_{\infty}$ moves continuously on $(P_{\infty})^c$.
This follows from the fact that
$\Phi$ is well defined and continuous in both $z$ and $a$ in
$\hat{\C} \sm \D \times (P_{\infty})^c$, analytic in $\hat{\C} \sm \D
\times (\oli{P}_{\infty})^c$. It follows that  $\oli{A}_{\infty}$
 moves continuously in the Hausdorff sense for $a \in (P_{\infty})^c$.
}


\end{proof}

\begin{cor}
\label{bubble lc}
The statement of \lemref{Ra} holds for $a\in (\bar P_\infty)^c$. Moreover,
for every such $a$, the bubbles of $R_a$ have locally connected boundaries.

\end{cor}
\begin{proof}
Consider a mapping $R_a$ with the parameter $a\in(\bar P_\infty)^c$ contained in
a capture component. Since $R_a$ is a hyperbolic mapping, the boundary of
every $A_\infty^a$ is locally connected by the standard considerations.
The second claim follows. The first claim is now immediate.
\end{proof}

\subsection{Internal parameter rays.}

If $P$ is a capture component of generation $n \geq 1$, for $t \in P$ let
$g_n(t)=\Phi_{t}(R^n(-1,t)))$, so that $g_n$ maps $a \in P$
to the B\"ottcher coordinate for $R_a^n(-1)$ in $A_{\infty}$. The function
$\xi_n$ a rational function and has a pole of finite order at the center of every capture component (later we show that it is in fact a simple pole). We proceed with the following definition.

\begin{defn}
An {\em internal parameter ray of angle $\th$} is a connected component of the set
\[
\{ g_n^{-1} (re^{2 \pi i \th}): r > 1 \}.
\]
\end{defn}

\begin{Lem} \label{land1}
Let $P$ be a parabubble with  $\gen(P)=n \geq 2$, and
let $\th\in\T$ be periodic (pre-periodic) under
doubling.
Then an internal parameter ray of $P$ with angle $\th$ lands at a point
$a_0\in\partial P$. Moreover, the point
$$p(a_0)=R_{a_0}^n(-1)$$ is a
repelling periodic (pre-periodic) point on the boundary of $A_{\infty}$.
\end{Lem}
\begin{proof}
To fix the ideas, we assume that $\th=0$ so that $p(a_0)$ is the
repelling fixed point where $A_{\infty}$ and $A_0$ touch.
Set $$\ga_n(t)=g_n^{-1}(te^{2\pi i\th}),\text{ for }t>1,$$
where we assume that $g_n^{-1}(te^{2\pi i\th})$ belongs to a chosen connected component of $\{ g_n^{-1} (re^{2 \pi i \th}): r > 1 \}$.
We want to show that $\lim\limits_{r \raw 1^+} \ga_n(r)$ exists and is
equal to $a_0$. First note that
\begin{equation} \label{deineq}
|\Phi_a^{-1}(r) - p(a)| \leq \de(r),
\end{equation}
where $\de(r) \raw 0$ as $r \raw 1$, which follows by Lemma \ref{rat}.
Also, note that the left hand side of (\ref{deineq}) is a continuous
function of both $a$ and $r$ on $\oli{P} \times (1,\infty)$.
This implies that $\Phi_a^{-1}(r) \raw p(a)$ uniformly as $r \raw 1$
on $\oli{P}$.

Therefore, for $a=\ga_n(r)$,
\begin{equation} \label{deineq2}
|\Phi_{\ga_n(r)}^{-1}(r) - p(\ga_n(r))| \leq \de(r),
\end{equation}
where $\de(r) \raw 0$ as $r \raw 1$. Now, for $|a-a_0| \leq \varep$, we have
$|R_a^n(-1)-p(a)| \leq \varep'(\varep) \raw 0$, as $\varep \raw
0$. On the other hand, since the zeros of $|R_a^n(-1)-p(a)|$ are
isolated, we can find a $C > 0$ such that if $0 < \varep \leq |a-a_0|
\leq C$, then $|R_a^n(-1)-p(a)| \geq \varep'$.

If $\ga_n(r)$ does not land at $a_0$, take an
$a \in \ga_n(r) \sm B(a_0,\varep)$, where $r$ is
sufficiently close to $1$, so that (\ref{deineq2}) holds
for $\de(r) \leq \varep'/2$. But since $|a-a_0| \geq \varep$ we have
$|R_a^n(-1)-p(a)| \geq \varep'$, for $a=\ga_n(r)$, which is a
contradiction. Hence $\ga_n(t)$ must land at $a_0$.
\end{proof}

The landing property for periodic parameter rays in $P_{\infty}$
follows from the standard theory in e.g. \cite{CG}, Theorem 5.2:
\begin{Prop}
If $\th$ is rational then the internal parameter ray of angle $\th$ in
$P_{\infty}$ lands at a parameter $a \in \partial P_{\infty}$. Moreover, if $\th \neq 0$ is periodic then $R_a$ has a parabolic cycle and if $\th$ is strictly preperiodic then $R_a$ is a postcritically finite map.
\end{Prop}

Consider the conjugacy $\phi$ from \lemref{capture conj}. We have the following:

\begin{Lem} \label{boundary}
Let $P$ be a parabubble of generation $n \geq 2$ and address
$\si$.
\begin{itemize}
\item[(I)] There exists a unique bubble $W\in K_\bas$ such that the following holds.
Let $a\in P$ and denote $B_a$ the bubble of $R_a$ which contains the critical value $-a$.
Then $\phi^{-1}(B_a)=W$.
\item[(II)] On the other hand, for each bubble $W \in \RR$, there exists a unique parabubble $P$ such that for any $a \in P$ we have $\phi(B_a)=W$, where $-a \in B_a$.
\item[(III)] Moreover, $a\in \partial P$ if and only if $-a \in \partial B_a$.
\item[(IV)] The parabubble $P$ is an open set, has a unique center, and is simply connected.
\end{itemize}
\end{Lem}

\begin{proof}
The first and third claim are immediate consequences of \lemref{above}. The same lemma implies that $P$ is an open set.
%

We have $\xi_n(a)=R^n(-1,a) \raw \partial A_{\infty}$ as $a \raw
\partial P$ by Lemma \ref{boundary}, so $\Phi_a \circ \xi_n \raw
\partial\D$ as $a \raw \partial P$. By the Argument Principle, this means that
every capture component $P$ is mapped by $\Phi_a \circ \xi_n$ onto
$\hat{\C} \sm \D$ as a $d\longrightarrow 1$ covering. We want to show that $d=1$.


Let $P$ be a parabubble of generation $n$, and $F$ the corresponding bubble for $R_a$ in which $-a$ lies. Note that the map $\phi$ in Lemma \ref{capture conj} is an injection of all bubbles of generation $\leq n$. Hence we can define $B=\phi^{-1}(F)$.
The root of $B$ then is a landing point $x$ of an internal ray of $B$
with angle $\th=0$ (by Lemma \ref{rat}). The
predecessor $C$ touches $B$ at $x$. It follows from Lemma \ref{land1}
that an
internal parameter ray with angle $\th=0$ in $P$ will land at a
parameter $a$ such that $R_a^n(-1)$ is the unique repelling fixed point on the boundary of $A_{\infty}$. It follows that there is a corresponding parabubble $Q$ to $C$ (in the same way as $P$ corresponds to $B$), such that $P$ touches $Q$ at $a$. Moreover, $gen(Q) < gen(P)$, since $gen(C) <gen(B)$. Proceeding in this way
we see that for every parabubble $P$, there is
a finite sequence of internal parameter rays connecting the center of $P$ with a point on $\partial P_{\infty}$.

Reversing this process we also see that for every bubble $B$ in the right basilica $\RR$ there is a corresponding parabubble $P$, in the sense that if $F$ is the bubble for $R_a$ in which $-a$ lies, then $B = \phi^{-1}(F)$. We cannot have such correspondence to the left basilica simply because $a=0$ is a singularity for the family $R_a$ and no sequence of parabubble rays can end there.

We have to prove that there is one and only one bubble in the right basilica corresponding to every parabubble. By Lemma \ref{xi-count} the only thing we have to show is that it is impossible to have one parabubble $P$ corresponding to two different bubbles $B_1$ and $B_2$ in the right basilica. This would imply that the parabubble has two distinct centers. By the $\la$-lemma of \cite{MSS},
any two centers in the same parabubble $P$ would correspond to
quasi-conformally conjugate rational maps. Since these maps would also
be postcritically finite, Thurston's Theorem implies
that a center is unique. Hence every parabubble corresponds to a unique bubble in the right basilica and (II) 
is proven.

Now, since the degree of $\xi_n$ coincides with the number of parabubbles of generation $n$, the Pigeonhole Principle implies that $\xi_n$ has a simple pole at the center of each parabubble of generation $n$.
By the Argument Principle,
$\Phi_a \circ \xi_n: P  \mapsto \hat{\C} \sm \D$ assumes every
value in $\hat{\C} \sm \D$ exactly once, so indeed $d=1$.
It follows that every capture component is simply connected.

\end{proof}
\noindent

By \lemref{boundary}, the mapping
$$\psi:a\mapsto -a\mapsto \phi^{-1}(-a)$$
is an injection from the capture locus of the family $R_a$ to
$\cir{\RR}$.
It is straighforward to extend this mapping to the roots of the
(para)bubbles, except for the roots contained in  the boundary of $P_{\infty}$.

Denote $\CC$ the union of capture components of the family $R_a$ and
$\CC_e = \CC \sm P_{\infty}$. Since dynamical bubbles may only touch at a single point, which is a
preimage of the fixed point where $A_{\infty}$ and $A_0$ as long as
$a \in \oli{P}_{\infty}^c$, our discussion implies:
\begin{Prop}
If $P$ and $Q$ are two parabubbles not equal to $P_{\infty}$, and
$P'=\psi(P)$, $Q'=\psi(Q)$, then the following holds:
\begin{enumerate}
\item $\oli{P} \cap \oli{Q} \cap (\oli{P}_{\infty})^c= \emptyset \Leftrightarrow \oli{P'} \cap
  \oli{Q'} = \emptyset$,
\item $\oli{P} \cap \oli{Q} \cap (\oli{P}_{\infty})^c$ is exactly one point $\Leftrightarrow$
$\oli{P'} \cap  \oli{Q'}$ is exactly one point,
\item $P = Q$ $\Leftrightarrow$ $P' = Q'$.
\end{enumerate}
Moreover,
$$\psi(\CC)=\cir{\RR}.$$
\end{Prop}
\noindent

Similarly to the notation for dynamical bubbles, if the intersection of
the closures of two parabubbles
 $$\oli{P} \cap \oli{Q} = \{a\}$$ is exactly one point and $\gen(P) > \gen(Q)$,
let us refer to $Q$ as the {\em predecessor} of
$P$ and $a$ as the {\em root} of $P$.

Let $\{a_j\}$ be the set of all touching points between parabubbles
not including those which lie on the boundary of $P_{\infty}$. The
above proposition implies that $\psi$ continuously extends to a
homeomorphism
$$\psi:\CC_e \cup \{ a_j\}\mapsto  \cir{\RR_e} \bigcup \(
\bigcup_{j=1}^{\infty} (f^{-j}(\al) \cap \RR_e) \).$$


\begin{defn}
Let $$\BB=\{F_k\}_0^\infty\subset K_\bas$$ be an infinite bubble ray
with angle $\angl(\BB)=\th\in (-1/3,1/3)$.
We call the corresponding sequence of capture components $\{P_k\}_0^\infty$, with
$\psi(P_k)=F_k$,  a
{\em parabubble ray} in $\CC$ with angle $\th$, and write $\angl (P)=\th$.

Similarly to the definition for dynamical bubble rays,
we define the {\em axis} for a parabubble ray $\PP$ to be the
union of the internal parameter rays $\ga_k,\ga_k' \subset P_k$ which land at
the points $\oli{P}_k \cap \oli{P}_{k-1} = x_k$ and $\oli{P}_{k+1}
\cap \oli{P}_k = x_k'$ respectively, starting from $\infty$.
\end{defn}

In the next section we show that certain infinite bubble
rays and parabubble rays land at a single point.

\section{Landing lemmas}

\subsection{Dynamical bubble rays.}
We begin with the following lemma.

\begin{Lem} \label{axis}
Assume that $\BB$ is a periodic infinite bubble ray
$\BB$ such that the axis is disjoint from the
closure of the postcritical set. Then the axis $\ga$ for $\BB$ lands
at a single periodic point which is either repelling or parabolic.
\end{Lem}
\begin{proof}
Let $\La$ be the closure of the postcritical set and let $S$ be the
set of cluster points for $\ga$. If the period of $\ga$ to itself is
$n$ then $R^n$ maps $\La \cup S$ into itself. Hence $R^{-n}$
can be lifted by the universal covering $\D$ of $\hat{\C} \sm (\La \cup
S)$ to a map $\hat{f}: \D \mapsto \D$ such that $\hat{f} (\D) \subset \D$
is a strict inclusion. Hence $R^n$ is strictly expanding with respect
to the Poincar\'e metric on $\hat{\C} \sm (\La \cup S)$.

Since $\ga$ is invariant under $f$ we can take a starting point $x_0
\in \ga$ and set $f(x_0)=x_1$, and $x_k=f(x_{k-1})$. Let $\ga_k$ be
the part of $\ga$ between $x_k$ and $x_{k+1}$. The
hyperbolic distance between $x_k$ and $x_{k+1}$ descreases as $k$
increases. Take a point $p \in S$. Then the hyperbolic distance from
any point on $\ga$ to $p$ is infinite, since $S$ is contained in the
boundary of the hyperbolic set $\hat{\C}\sm (\La \cup S)$.
Since the hyperbolic length of $\ga_k$ decreases for
increasing $k$, any neighbourhood $N$ of $p$ has the property that
there is a smaller neighbourbooh $N' \subset N$ such that if $\ga_k
\cap N' \neq \emptyset$ then $\ga_k \subset N$. But this means that
$f(N) \cap N \neq \emptyset$. So $p$ has to be a fixed point. Since
$S$ is connected, $S$ must contain only this point.
By the Snail Lemma, $p$ must be a parabolic or repelling point
(cf. \cite{Milnor1}, Lemma 16.2).
\end{proof}

We next prove that the axis of a periodic (or preperiodic) bubble ray cannot
accumulate on some bubble.
\begin{Lem} \label{noacc}
Let $\BB$ be a periodic infinite bubble ray for which the axis is disjoint
from the closure of the postcritical set. Then the axis of $\BB$
cannot accumulate at some bubble.
\end{Lem}

\begin{proof}
Without loss of generality, in order to reach a contradiction, it
suffices to suppose that the axis of $\BB$ accumulates at $\partial
A_{\infty}$. Since
$\BB$ is periodic, we know from Lemma \ref{axis} that
the axis for the bubble ray $\BB$ lands at a single periodic point
$p$ on the boundary of $A_{\infty}$. The bubble ray $\BB$
then encloses a domain $D$ whose boundary is a connected part $I \neq
A_{\infty}$ of $\partial A_{\infty}$ and half of the boundary of all
the other bubbles in $\BB$.
Since $p$ and $\BB$ is fixed under some iterate $n$ we have
that $D$ is invariant under $R^n$. This means that any bubble $B$ in
$D$ must never be mapped into $A_0 \cup A_{\infty}$, since this set
lies outside $D$ (the fact that there exists some bubble in $D$ is
obvious). This is clearly impossible, since bubbles by
definition are preimages of $A_{\infty}$.
\end{proof}

We are now in position to prove a landing lemma for periodic or
preperiodic bubble rays.

\begin{Lem} \label{landingd1}
Assume that $\BB$ is periodic infinite bubble ray, for which there
exist an $N$ such that all
bubbles in $\BB$ of generation at least $N$ are disjoint from
the closure of the postcitical set. Then $\BB$ lands at a single point.
\end{Lem}

\begin{proof}
Assume that $\BB$ is periodic of period $q$.
We have seen (Lemma \ref{axis} and Lemma \ref{noacc}) that the axis
$\ga$ of the bubble ray must land on a periodic point $x$.




Since the postcritical set $\La$ is disjoint from any bubble $B$ in $\BB$
with $Gen(B) \geq N$,
we have an annulus $R$ around this $B$ of some definite
modulus $m > 0$ such that there are well defined inverse branches
of $R_a^{-q}$ on $R \cup B$, where $R_a^{-qn}(B) \in \BB$ for all $n
\geq 0$.
This means that the lengths of the
$\ga_k$ in the proof of Lemma \ref{axis} are commensurable with the
diameter of the corresponding bubbles $F_k$, by the Koebe Distortion
Lemma.
Hence the bubble ray $\BB$ converges to the same periodic point
as the axis $\ga$ lands on.

\end{proof}


\subsection{Orbit portraits for $R_a$.}
We have seen in Section \ref{raysandaxes} that bubble rays have angles
inherited from the angles of external rays in the basilica (although
these angles are not always well defined, as in the capture case for
instance). With the theory about orbit portraits for quadratic
polynomials in Section \ref{orbit}
in mind, it is now straightforward to define an orbit portrait for
$R_a$.

\begin{defn}
Let $x_1,x_2,\ldots,x_p$ be a (repelling or parabolic) periodic orbit, where
$R_a(x_i)=x_{i+1}$, $R_a(x_p)=x_1$. Assume that there
are a finite number of periodic infinite bubble rays landing on
$x_i$, with well defined angles; Let $A_i$ be the
corresponding angles for the bubble rays landing at $x_i$.
Then {\em the orbit portrait for $R_a$} is the set
$\OO=\{A_1,A_2,\ldots,A_p\}$.
\end{defn}


Given two angles $\th_1\neq\th_2$ we let $\circarc{\th_1}{\th_2}\subset \T$ be
the arc of the unit circle
swept by going in counter-clockwise direction from $\th_1$ to $\th_2$. We
say that $\th$ lies between $\th_1$ and $\th_2$ if $\th\in \circarc{\th_1}{\th_2}$.



Before we state the next lemma we make some more definitions.


\begin{defn} \label{innerbound}
Let $\BB_1$ and $\BB_2$ be two bubble rays starting from $A_\infty$
with well defined angles $\th_1$ and
$\th_2$ and axes $\ga_1$ and $\ga_2$. Assume that $\BB_1$ and $\BB_2$
land at a common point $p$.
Denote $D$ the domain bounded by the axes $\ga_1$, $\ga_2$ which does not contain any
bubble rays with angles in $\T\setminus [\th_1\circlearrowleft\th_2]$.

Define the {\em outer boundary} of the sector bounded by $\BB_1$, $\BB_2$
as the union of the arcs of the boundaries of the
bubbles in these two bubble rays lying outside $D$ together with their
endpoints. Similarily, define the
{\em inner boundary}.
We say that $z\in\hat\C$
lies {\em between} $\BB_1$ and $\BB_2$ if $z \in D$ and $z \notin
\oli{\BB}_1 \cup \oli{\BB}_2$.
\end{defn}

This notion of being {\em between} two bubble rays also
makes sense for bubble rays even if $\BB_1$ and $\BB_2$ do not land on a common
point.
\begin{defn}
Assume that the two bubble rays $\BB_1$ and $\BB_2$ have
intrinsic addresses $s(\BB_1)=(x_0,x_1,\ldots,)$ and
$s(\BB_2)=(y_0,y_1,\ldots)$ respectively.
We say that an infinite bubble ray $\BB$, with intrinsic address
$s(\BB)=(z_0,z_1,\ldots)$,
lies {\em between} $\BB_1$ and
$\BB_2$ if $y_i \leq z_i \leq x_i$ for all $i \geq 0$.
Equivalently,
the angle $$\angl(\BB)\in [\angl(\BB_1)\circlearrowleft \angl(\BB_1)].$$
\end{defn}

This definition also makes sense for parabubble rays in an exactly
analoguous way.

\begin{Lem} \label{landingd2}
Let $\OO=\{A_1,\ldots,A_p\}$ be a formal orbit
portrait with $v_\OO\geq 2$ and let $I = \circarc{t_{-}}{t_{+}}$ be its characteristic arc.
If the formal orbit portrait $\OO$ is realisable
by some $R_a$ then $-a$ cannot lie on the outer boundary of a bubble ray
with angle $t_{-}$ or $t_{+}$.
\end{Lem}

\begin{proof}
Since $a \in \mat$, in the case when  $-a$ belongs to the boundary of some bubble,
we have a conjugacy $\phi$ from \propref{marking-mated} between the dynamics of
$f_\bas$ on the interior of $K_\bas$ and that of $R_a$ on its Fatou set.
Now suppose $\OO$ is realised and let
$A_i$ be the set of angles of the bubble rays landing at $x_i$.

Assume that $A_2$ contains the characteristic arc.
Let $\BB_{-}$ and $\BB_{+}$ be the bubble rays corresponding to the angles
$t_{-},t_{+} \in A_2$ and let $\AA_{+},\AA_{-}$ be the bubble rays
corresponding to the critical arc in $A_1$, i.e. so that $\AA_{-}$ and
$\AA_{+}$ are mapped onto $\BB_{-}$ and $\BB_{+}$ respectively. Also,
let $D$ be the domain enclosed by the axes of $\BB_{-}$ and
$\BB_{+}$.

There are two more preimages of $\BB_{-}$ and $\BB_{+}$, call them
$\AA_{-}',\AA_{+}'$ respectively. Also, let $a_{-}=\angl(\AA_{-}),
a_{+}=\angl(\AA_{+})$ and $a_{-}'=\angl(\AA_{-}'), a_{+}'=\angl(\AA_{+}')$.
Since $I_c=(a_{-},a_{+})$ is the critical arc in $A_1$ we have that both
$a_{-}',a_{+}'$ lies entirely inside $I_c$, and thus the bubble rays
$\AA_{-}', \AA_{+}'$ lies entirely inside the domain $D_c$ enclosed by the
axes for the bubble rays forming the critical arc.

Now, assume that $-a$ lies on an outer boundary of a bubble in
$\BB_{+}$ (the proof is the same if $-a \in \BB_{-}$). Then the
critical point $-1$ must belong to
a bubble in $\AA_{+}$. Since $R_a$ is $2-1$ in a neighbourhood of $-1$
and orientation preserving, we have that the bubble ray $\AA_{+}'$
must touch $\AA_{+}$ at $-1$. Since $-a$ is outside $D$, this implies
that $-1$ must be outside $D_c$.
Thus $\AA_{+}'$ must be outside $D_c$, which is a contradiction.
\end{proof}

The following lemma tells us when a specific orbit portrait is
realised.

\begin{Lem}[{\bf Realization of orbit portraits}] \label{realise}
Let $\OO=\{A_1,\ldots,A_p\}$ be a formal orbit portrait with
a characteristic arc $\II = \circarc{t_{-}}{t_{+}}$. Let $P_{t_{-}}, P_{t_{+}}$
be the corresponding parabubble bubble rays, with angles $t_{-}$ and
$t_{+}$ and assume that $a$ belongs to a parabubble $P$ between $P_{t_{-}}$ and
$P_{t_{+}}$. Then the orbit portrait $\OO$ is realised by $R_a$.
\end{Lem}
The proof follows that of Lemma 2.9 in \cite{Milnor2}.
\begin{proof}
Note that all infinite bubble rays with angles in any $A_j$ are
well defined since their forward images do not intersect the critical value.

Let $\La$ be the closure of the postcritical set for $R_a$ and let
$\rho(z)$ be the induced hyperbolic metric on $\hat{\C} \sm \La$.
Let $C$ be the critical bubble containing $-1$ and $V$ the critical
value bubble containing $-a$. There is a unique finite bubble ray
ending at $V$. Its preimage is two finite bubble rays $\BB_1$ and
$\BB_2$ both ending at $C$. Their axes $\ga_1$ and $\ga_2$ join in $C$
and form a closed simple curve in $\hat{\C}$.

Take a hyperbolic disk $D \Subset C$ which covers the critical point
$-1$ and let
\[
L = \bigcup_{k=0}^{\infty} R_a^k (\ga_1 \cup \ga_2 \cup D)
\]
It is easy to see that the complement of $L$ is two topological disks
$U_1$ and $U_2$.

We have $R(L) \subset L$ and $\La \subset L$. Moreover
$\dist(\La,L^c) \geq \varep > 0$ for some definite $\varep
> 0$. It is easy to check that the $n$th preimages of $U_1$ and $U_2$
consist of $2^{n+1}$ topological disks.

Moreover, all preimages of the $U_j$ will be on a definite
distance $\varep >0$ from $\La$ so we have a uniform constant
$c=c(\varep) > 1$ so that
\[
\rho(R(x),R(y)) \geq c \rho(x,y)
\]
for $x,y$ lying in any of these preimages of $U_i$.
It follows that the preimages
of $U_i$ shrink to points. Thus the symbol sequence of some point with respect
to the initial partition $L$ is unique. In particular the landing
points of the periodic bubble rays in $\OO$ will have the same symbol
sequence if and only if they land at a common point.

To show that $\OO$ is indeed realised it now suffices to show that all
the landing points of the bubble rays with angles in $A_j$ lie
entirely in one of the components $U_i$. Since they are
mapped onto each other they will have the same symbol sequence in that
case.

The preimages of
the characteristic arc $\circarc{t_{-1}}{t_{+}}$ under the doubling map
will be two smaller arcs $I'$ and $I''$ at the end of the critical arc.
Since every $A_j \in \OO$ cannot have any element in $I'$ or $I''$ we
have that all bubble rays corresponding to angles in $A_j$ are
completely contained in $U_1$ or completely contained in $U_2$. Thus
all the angles in every $A_j$ have the same symbol sequence, so they
land on a common point, and so $\OO$ is a realised bubble portrait.
\end{proof}

\begin{Lem} \label{highlander}
Assume that $R_a$ has a parabolic fixed point $z_0$, with
$R_a'(z_0)=e^{2 \pi i p/q}$, where $p/q \in \Q$ with $(p,q)=1$.
Then there are precisely $q$ periodic bubble rays $\BB_j$, $j=1,\ldots,q$,
landing at $z_0$. These bubble rays are mapped onto each other
under the action of $R_a$, with combinatorial rotation number $p/q$.
\end{Lem}
\begin{proof}
For simplicity, consider the mapping $R_a$ with $a=32/27$ which has a simple
parabolic with eigenvalue $1$. After a suitable change of coordinates shifting the
fixed point to the origin, this mapping takes the form
$$\zeta\mapsto\zeta+\zeta^2+\OO(\zeta^3)$$
in a neighborhood of $\zeta=0$.
Denote $\AA$ and $\RR$ the attracting and repelling petals of $R_a$ correspondingly.
Note that Montel's Theorem guarantees that the repelling petal contains a bubble $B$.

Now, $B$ is the end of some finite bubble ray $\CC_F$.
Taking the preimages of the bubble ray $\CC_F$ we get a sequence of
bubble rays $\CC_k=R_a^{-k}(\CC_F)$, whose ends will converge to
$z_0$.

Since preimages will increase the generation and since
there are finitely many finite bubble rays of any fixed generation,
for any $N$ there must be some bubble $\BB_0$ of generation $N$ contained in
infinitely many $\CC_k$. Let $\CC_{k_0} \supset \BB_0$ be the $\CC_k$
containing $\BB_0$ with lowest generation and $\CC_{k_1} \supset \BB_0$
the second lowest. Then $R_a^m(\CC_{k_1}) = \CC_{k_0}$ for some $m
\geq 1$ and the preimage of $\BB_0$ under $R_a^m$ is a longer bubble
ray $\BB_1 \supset \BB_0$. Moreover, $\BB_1 \subset \CC_{k_1}$. Taking
further preimages of $\BB_1$ under the same branch $f=R_a^{-m}$ we
get a sequence $\BB_n$ of nested finite bubble rays such that $\BB_n
\subset \CC_{k_n}$. Moreover, the ``difference'' between $\BB_n$ and
$\CC_{k_n}$, i.e. the number of bubbles in $\DD_n=\CC_{k_n} \sm \BB_n$ is a
fixed constant $K$ for all $n$. The bubble $\DD_n$ is also a preimage
of the starting set $\DD_0$ under $f$. Since the postcritical set $\La$
accumulates on $z_0$, it is disjoint from
$\DD_n$. Thus there is a neighbourhood around all bubbles in $\DD_0$
where $f^n$ is defined for all $n \geq 0$.
Now, the Koebe Distortion Lemma implies that all bubbles in
$\DD_n$ shrinks to points, namely the parabolic fixed point
$z_0$, since one of them, namely the end of $\CC_{k_n} \supset \DD_n$,
converges to $z_0$. Hence there is a subsequence of bubbles in $\BB_n$ which
converge to $z_0$ (but we do not know a priori that the bubble ray
itself will converge to $z_0$).

However, by construction, the bubble ray $\BB=\cup_n \BB_n$ is periodic.
We can now apply Lemma \ref{landingd1} to $\BB$, which shows that
$\BB$ lands at a single point, which must be equal to $z_0$.

Let us show that the period of $\BB$ is $1$. A priori, $\CC$ is
periodic with a period which divides $m$. Assume the period is $p \neq
1$ and that $R_a(\BB_j)=\BB_{j+1}$ for $1 \leq j \leq p-1$,
$R_a(\BB_p)=\BB_1$.
By simple combinatorial considerations (see e.g. \cite{Milnor2}),
these bubbles form their own orbit portrait. But this means that some
point $z \in \RR \sm \AA$ in the domain bounded by
two consecutive bubble rays $\BB_j$ and $\BB_{j+1}$, will be mapped
into $\AA$, which is impossible. Hence $p=1$.
\end{proof}

\subsection{Parameter bubble rays.}

Let us first note the following evident statement:
\begin{lem}
\label{orbit persist}
Assume that an orbit portrait $\OO$ is realized for some rational map $R_a$
by bubble rays landing at a repelling orbit $\{x_i\}$. Let $a_t$, $t\in[0,1]$
be a continuous path with $a_0=a$ along which the corresponding periodic orbit $\{x_i^t\}$
remains repelling.
Assume further that for every $t$ no iterate of the critical
value $-a_t$ is contained in the boundary of a bubble ray with angle $\gamma\in\OO$. Then
the orbit portrait $\OO$ is realized for all $R_{a_t}$.

\end{lem}


The following Proposition has an analogue in \cite{Milnor1},
Theorem 4.1 (and Lemma 4.2). Since the proof is completely similar, we
omit it.
\begin{Prop}[Milnor; Parameter Path] \label{mil}
Given a parameter $a_0$ such that $R_{a_0}$ has a parabolic fixed
point $z_0$ with combinatorial rotation number $p/q$ and
an orbit portrait $\OO$ (from Lemma \ref{highlander}).
Then there is a path $\ga$ emerging from $a_0$ in parameter space so
that $a \in \ga$ implies that $R_a$ has a repelling fixed point
$z=z(a)$ with orbit portrait $\OO$ and an attracting periodic orbit
with period $q$, close to $z(a)$.
\end{Prop}

The set $A$ of parameters where the attracting periodic orbit in the above
lemma exists, is bounded by a finite number of analytic curves.
Indeed,
\[
A = \{ a : |(R^q)'(z_i(a),a)| < 1 \}.
\]
The condition $|(R^q)'(z_i(a),a)| = 1$ represents an analytic curve
with a finite number of singularities.
We conclude that there
is a ``wedge'' $\ti{W}$, that is, two analytic curves $\ga_1$ and $\ga_2$ which
meet at $a_0$ such that for a small neighbourhood $B(a_0,\varep)$, an
open set $E$ bounded by $\ga_1$, $\ga_2$ and $\partial B(a_0,\varep)$ has
the property that inside $E$, we have $\OO$ realised and
$z_i(a)$ is an attracting periodic orbit of period $q$ (as in the
above lemma).

By \lemref{orbit persist} and \lemref{landingd2}
the parabubble rays $\PP_{t^+}, \PP_{t^-}$ lie outside of the wedge
$\ti{W}$.

\begin{Lem} \label{accum}
Let $a_0$ be as in the above lemma and assume that $(t^+,t^-)$ is the
characteristic arc for $\OO$. Then for any $\varep > 0$, we have
$B(a_0,\varep) \cap \PP_{t} \neq \emptyset$, for at least one
$t=t^+, t^-$, where
$\PP_{t^+},\PP_{t^-}$ denote the parabubble rays with angles $t^+,t^-$
respectively.
\end{Lem}
\begin{proof}
Assume the contrary. Then there is an $\varep > 0$ such that
$B(a_0,\varep)$ is disjoint from the parabubble rays $\PP_t$ for
$t=t^+,t=t^-$. By the above argument, and Theorem \ref{mil},
the orbit portrait $\OO$ is realised in $B(a_0,\varep) \cap N$,
where $N = \{a: |R_a'(\al(a))| > 1 \}$, and $\al(a)$ is the (local)
continuation of the parabolic fixed point
$z_0$ (this is possible if the multiplier is $\neq 1$).
Hence there is a parameter $a_1 \in
B(a_0,\varep)$, such that $R_{a_1}$ also has a parabolic fixed
point $z_1$.

But since the combinatorial rotation number is changed for $a_1$
the new wedge $\ti{W}_1$ emerging from $a_1$ has to exhibit a
different orbit portrait $\OO_1$. But $\ti{W}_1$ must intersect
$B(a_0,\varep) \cap N$, and so both orbit portraits
$\OO_1$ and $\OO$ are realised, which is impossible. The lemma follows.
\end{proof}

\begin{Prop} [{\bf Parabubble wakes I}] \label{paralanding}

Let $a_0$ be such that $R_{a_0}$ has a parabolic fixed point
$z_0$ with eigenvalue $R_{a_0}'(z_0)=e^{2\pi ip/q}$, $(p,q)=1$.
Denote $\OO=\{\{\th_1,\ldots,\th_q\}\}$ the orbit portrait from
Lemma \ref{highlander}, and let $\II=\circarc{t_{-}}{t_{+}}$ be
its characteristic arc. Then the corresponding parabubble rays with angles
$t_{+}$ and $t_{-}$ land on $a_0$.
\end{Prop}

\begin{proof}
The standard considerations of parabolic dynamics imply that
$$R_{a_0}^q(z)=R(z)=(z-z_0)+b(z-z_0)^{q+1}+\OO((z-z_0)^{q+2}),$$
 for some $b \neq 0$.
For $a$ close to $a_0$ the fixed point $z_0$ will bifurcate into
$q+1$ fixed points (for $R_a^q$) $z_k(a)$, which are analytic in a
neighbourhood of
$a_0$, and where $z_k(a_0)=z_0$, for $k=1,\ldots,q+1$. One of these
fixed points must be a fixed point for $R_a$ as well if $q \geq 2$,
while the other fixed points (for $R_a^q$) are all repelling, indifferent or
attracting. By Lemma \ref{accum} there must be a subsequence of
parabubbles $P_{n_k} \subset \PP_{t^+}$ (or $\PP_{t^-}$) such that
$P_{n_k} \cap B(a_0,\varep) \neq \emptyset$, for all $k \geq
N(\varep)$. Hence, for sufficiently large $k$, if $a_1 \in P_{n_k}$,
then $-a_1 \in B^{n_k}$, where $B^{n_k}$ is the corresponding
dynamical bubble in the bubble ray $\BB_{t^+}$, i.e. with same address as
$P^{n_k}$.
Since $a_1$ is a capture parameter the fixed points $z_i(a_1)$ (under
$R_a^q$) cannot be attracting. They cannot be neutral so they must be
repelling.

We now use the standard theory of parabolic bifurcation (see for ex
\cite{Shishikura1} Section 7, \cite{Shishikura2},
\cite{Orsay-notes}). For a suitable small perturbation, we get $q$ fundamental
domains $S_{+,a}^k$ and $S_{-,a}^k$, $1 \leq k \leq q$, for the
repelling and attracting petals respectively for the perturbed map
$R_{a}$. They have the property that
$$S_{+,a}^k \cap S_{-,a}^k = \{\al(a),z_k(a) \}.$$
Moreover, there exist analytic functions $\Phi_{+,a}^k,
\Phi_{-,a}^k$ (the perturbed Fatou coordinates)
which are defined and injective in a neighbourhood of
$\ti{S}_{+,a}^k=S_{+,a}^k \sm \{\al(a), z_k(a) \}$ and $\ti{S}_{-,a}^k=S_{-,a}^k \sm
\{\al(a), z_k(a) \}$ respectively, and conjugate the dynamics
of $R_a^q$ to that of the unit translation. With a choice of
normalization, these coordinates will vary locally analytically
with $a$.

 If $z \in \ti{S}_{-,a}^k$, then there is
an $n \geq 1$ such that $R_a^{qn}(z) \in \ti{S}_{+,a}^k$, and for the
smallest such $n$,
\[
\Phi_{+,a}^k(R_a^{qn}(z)) = \Phi_{-,a}^k(z) - \frac{1}{\be(a)} + n+\operatorname{const},
\]
where $\be(a)=\be_k(a)$ is an analytic function in a punctured neighbourhood of
$a_0$, $\be(0)=0$, defined by
$$(R_a^q)'(\al(a)) = e^{2 \pi i \be(a)}.$$

Denote $C_+^k$, $C_-^k$ the {\'E}calle-Voronin cylinders, obtained as the
quotients
$$C_+^k=S_{+,a}^k\mod R^q_a\simeq\C/\Z,\;C_-^k=S_{-,a}^k\mod R^q_a\simeq\C/\Z.$$
We get that for $z \in C_{+}^k$,
\begin{equation} \label{eqq}
\Phi_{-,a}^k \circ R_a^{qn} \circ (\Phi_{+,a}^k)^{-1}(z) = z +
\frac{1}{\be(a)} \mod \Z.
\end{equation}
The function $$\tau_a(z)=\frac{1}{\be(a)} + z \mod \Z$$ viewed as an isomorphism
$C_{+}^k\mapsto C_{-}^k$ is called the transit map.

Now let us fix some $k$ so that the critical point $-1$ belongs the
the $k$th attracting petal. For simplicity let us drop the indices $k$
in the above discussion and only focus on these particular Fatou
coodinates. Then for any prescribed bubble $B_l$ in $\BB_{t^+}$ we can
find a parameter $a \in B(a_0,\varep)$ such that $R_a^{qn}(-1) \in
B_l$, for some $n \geq 1$, $n=n(l)$.

Fix $a=a_1$ as above. For this specific perturbation, we already have $-a \in
B^{n_k}$, so we know that $n=1$ in (\ref{eqq}). The bubbles
$B_n$ move holomorphically and with uniformly bounded distortion in
the Fatou coordinates for $a$ in some disk $B(a_0,\varep)$
(the lifted dynamical bubbles in $\BB_{t^+}$, in the Fatou
coordinates, are all unit translates of each other). The function
$\tau_a(-1)=1/\be(a) + z \mod \Z$ has derivative
\[
\partial_a \tau_a(z) \sim D \frac{1}{(a-a_0)^{m}} = \frac{-m}{(a-a_0)^{m+1}},
\]
for some $m \in \Q$, $m > 0$.
This, and the distortion considerations, imply that  $P_{n_k}$ converge to $a_0$.

It remains to show that {\em all} parabubbles $P_l \subset \PP_{t^+}$
converge to $a_0$, instead of just a subsequence $l_k$. This follows
from the fact that $n=n(a)$ in (\ref{eqq}) is continuous function of $a$ which
only assumes integer values. Hence $n=1$ for all $l$ and $\PP_{t^+}$
lands on $a_0$. Of course a similar statement holds for $\PP_{t^-}$.


\end{proof}

Let us write $W=W(t^+,t^-)$ for the {\em parabubble wake} being
set of points between the parabubble rays from the above lemma.
Also, let $\OO=\OO(t^+,t^-)$,
be the corresponding orbit portrait. Note that the characteristic arcs
corrsponding to different orbit portraits around the fixed point are
disjoint.

\begin{figure}[ht]

\centerline{\includegraphics[height=0.42\textheight]{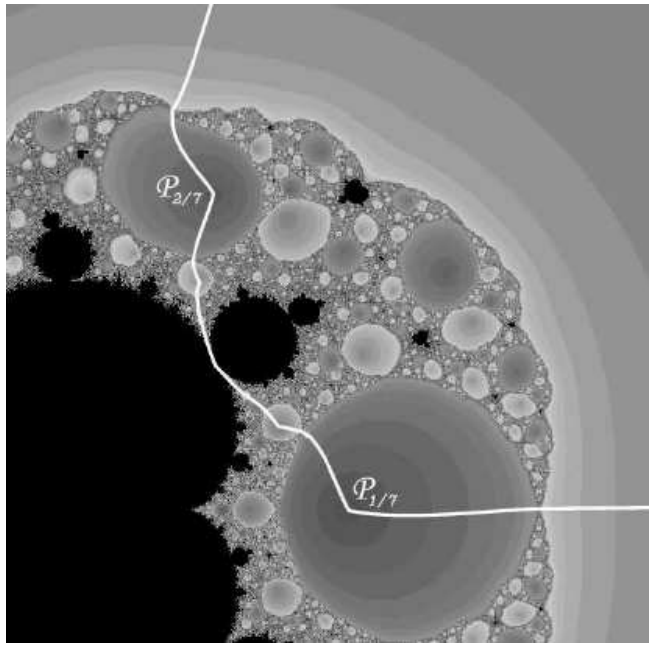}}
\caption{An example of a parameter wake, $W(1/7,2/7)$. The axes of the
parabubble rays $\PP_{1/7}$, $\PP_{2/7}$ which bound the wake are indicated.
Their common landing point is the parameter value $a_0$ for which $R_{a_0}$
has a parabolic fixed point $z_0$ with eigenvalue $e^{2\pi i/3}$. The orbit
portrait of $z_0$ is $\{1/7,2/7,4/7\}$.\label{wake figure}} 
\end{figure}

\begin{Lem}[{\bf Parabubble wakes II}]
\label{wakesII}
The parabubble rays in the above lemma cut out an open set in the complex
plane, called the {\em bubble wake} $W=W(t^+,t^-)$ such that $a \in W$ if
and only if $R_a$ exhibits the repelling orbit portrait $\OO=\OO(t^+,t^-)$.
\end{Lem}
\begin{proof}
By Lemma \ref{boundary} the set $A_{\infty}$ moves holomorphically and
the critical value $-a$ belongs to the boundary of a bubble if and
only if $a$ belongs to the boundary of the corresponding
parabubble.
By \lemref{orbit persist}
if for a single parameter $a \in W=W(t^+,t^-)$ the map $R_a$ realises the
orbit portrait $\OO=\OO(t^+,t^-)$, then the same is true for every parameter in $W$.

On the other hand, $\OO$ cannot be
realised for any parameter value outside $W$.
Indeed,  $\OO$ is not realised for $a$ in  any of the capture components outside
$W$, since this would imply that the critical value is outside the
characterisic arc.
\end{proof}

\section{A puzzle partition for $R_a$} \label{puzzle}

The idea of a puzzle partition for a Julia set originated in the work
of Branner and Hubbard \cite{BH}. It has been further developed by Yoccoz
(see e.g. \cite{Hubbard} and \cite{Milnor5}), to study the local connectedness
of the Mandelbrot set at Yoccoz parameters, and the local connectedness of the
corresponding Julia sets. We employ the Branner-Hubbard-Yoccoz approach to
maps of the family $R_a$ using partitions given by landing bubble rays.

\subsection{The Yoccoz puzzle for quadratic polynomials.}
Let us recall the main steps of Yoccoz' construction for a
quadratic polynomial $f_c$ without non-repelling orbits with a connected Julia set.
Let $\alpha$ stand for the dividing fixed point of $f_c$. It is the landing point
of a cycle of $q>1$ external rays of $f_c$.
Denote these rays $R_1,\ldots, R_q$.
Recall that the B{\"o}ttcher coordinate
\[
\Phi: \hat{\C} \sm K(f_c) \mapsto \hat{\C} \sm \D,
\]
conjugates $f_c$ to the dynamics of $z\mapsto z^2$.
Fix an arbitrary $r>1$ and let $E_r$ be the equipotential curve
$$E_r=\Phi^{-1}(\{ re^{2 \pi i \th}: \th \in [0,1] \}.$$
Let $U_0$ be the graph formed by
$$U_0=R_1\cup\cdots\cup R_q\cup E_r\cup\{\alpha\}.$$
The {\it puzzle pieces of depth $0$} are the bounded components of
$\C\setminus U_0$. Denote these $q$ topological disks $P_0^j,\;
j=0,\ldots q-1$. By definition, the Yoccoz' puzzle pieces of depth
$d\geq 1$ are the first preimages of the puzzle pieces of depth
$d-1$ under $f_c$.

\comm{
\begin{figure}
\psfrag{a}[][][0.7]{$\al$}
\psfrag{b}[][][0.7]{$-\al$}
\psfrag{c1}[][][0.7]{$Z_1^1$}
\psfrag{c2}[][][0.7]{$Z_2^1$}
\psfrag{d}[][][0.7]{$Y_0^1$}
\begin{center}
\includegraphics[scale=.5]{firstpuzzle.eps}
\end{center}
\caption{The Yoccoz puzzle at depth $1$, with $q=3$ external rays
  landing at $\al$.\label{puzzle1}} 
\end{figure}
}

\begin{figure}
\centerline{\includegraphics[width=1.1\textwidth]{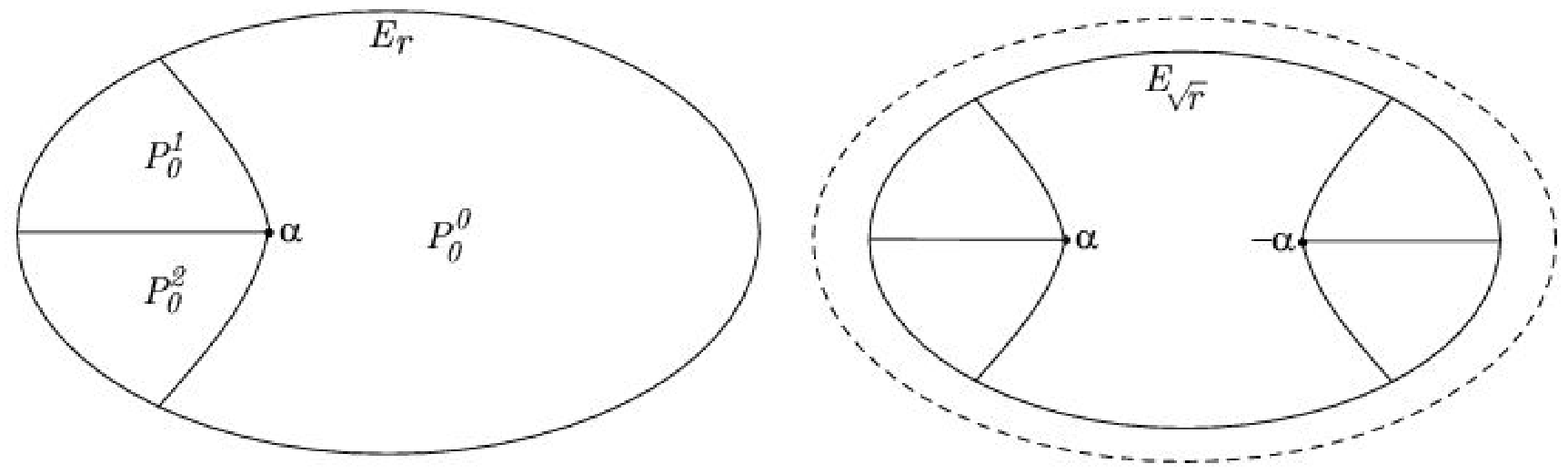}}
\caption{The Yoccoz puzzles of depths $0$ and $1$, with $q=3$ external rays
  landing at $\al$. \label{puzzle1}}
\end{figure}

What makes puzzle partitions of Julia sets so useful in the study of local
connectedness are the following two straightforward observations:

\begin{prop}
The following two properties hold:
\begin{itemize}
\item (Markov property)
any two puzzle pieces $P_d^j$ and $P_{d'}^{j'}$ are either disjoint,
or one of them is contained in the other;
\item the intersection $J_c\cap P_d^j$ is a connected set.
\end{itemize}
\end{prop}

\noindent
The Markov property allows us to make the following definition for any point $z\in J_c$
which is not a preimage of $\alpha$.
\begin{defn}
For any $z\in J_c$ with $\alpha\notin \cup f_c^n(z)$,
let $P_d(z)$ denote the puzzle piece of depth $d$ which contains $z$.
Let us also set
\[
A_d(z) = P_d(z) \sm \overline{P_{d+1}(z)}.
\]
\end{defn}

We will refer to $A_d(z)$ as an annulus, even though it may be degenerate.
The sequence of annuli $A_d(0)$ will be called the {\it critical annuli}.

The following is a consequence of Gr{\"o}tzch Inequality (see e.g. \cite{BH}):
\begin{lem}
\label{moduli diverge}
Let $A_i$, $i\in\N$  be a sequence of bounded conformal
annuli in the plane with simply-connected complementary components.
 Denote $W_i$ the bounded component of $\C\setminus \bar A_i$.
Assume that  $A_{i+1}\subset W_i$
and
$$\sum\mod A_i=\infty.$$
Then
$$\diam \left(\bigcap W_i\right)=0$$
\end{lem}

Yoccoz has demonstrated, in particular:

\begin{lem}
Assume that $f_c$ is non-renormalizable. Then
$$\sum\mod A_d(0)=\infty.$$
\end{lem}

\noindent
His proof uses the concept of a {\it tableau} developed by Branner and Hubbard \cite{BH}.
Below we extract a definition suitable for a generalization from \cite{Milnor5}.
To motivate some of the notation, fix a point $z\in J_c$, and consider its orbit under $f_c$:
\[
z=z_0 \mapsto z_1 \mapsto z_2 \mapsto \ldots.
\]
Note that the puzzle piece $P_d(z_j)$
is mapped onto $P_{d-1}(z_{j+1})$, either as a conformal isomorphism
or a branched double covering, depending on whether the piece $P_d(z_i)$
contains the critical point or not.

\begin{defn}
Let $S(z)$ be the largest integer $d \geq 0$, for which
$P_d(z)=P_d(0)$. If $P_d(z)=P_d(0)$ for all $d$, put $S(z)=\infty$, and if $P_d(z)
\neq P_d(0)$ for all $d$, put $S(z)=-1$.
\end{defn}

\noindent
We then distinguish the following three possibilities:

\begin{itemize}
\item
{\bf Critical case:} $d < S(z_i)$. Here the critical point lies in
$P_d(z_i)=P_d(0)$. Hence the annulus $A_d(z_i)$ is mapped onto its
image as am unbranched two-to-one covering. One easily deduces that
\[
\mod A_d(z_i) = \frac{1}{2} \mod A_{d-1}(z_{i+1}).
\]

\item
{\bf Off-critical case:} $d > S(z_i)$. Here the critical point is
outside $A_d(z_i)$ so that $A_d(z_i)$ is mapped conformally onto its
image $A_{d-1}(z_{i+1})$. Indeed,
\[
\mod A_d(z_i) = \mod A_{d-1}(z_{i+1}).
\]

\item
{\bf Semi-critical case:} $d=S(z_i)$. This means that the critical
point lies in the annulus $A_d(z_i)$, and
\[
\mod A_d(z_i) > \frac{1}{2} \mod A_{d-1}(z_{i+1}).
\]
\end{itemize}



\begin{defn}[{\bf A critical tableau}]
A critical tableau is a two-dimensional array of non-negative real numbers
$(\mu_{d,n})$, $d,n\geq 0$ together with a marking, formed according to a set of
rules given below.
Each position of the
tableau is marked as {\it critical}, {\it semi-critical}, or {\it off-critical}.
An {\it iterate} $\II$ in the tableau is a move in the north-western direction in the array:
$$\mu_{d,n}\underset{\II}{\longrightarrow} \mu_{d-1,n+1}.$$
The rules of a critical tableau are as follows.

\begin{itemize}
\item Every column of a tableau is either all critical; or all off-critical;  or has exactly one
semi-critical position $(d_0,n)$ and is critical above ($d>d_0$) and off-critical below.
The $0$-th column is all critical.

\item If $$\mu_{d,n}>0\text{ then }\II(\mu_{d,n})>0.$$
\noindent
\begin{tabbing}
Moreover, \= if $(d,n)$ is marked off-critical, then $\II(\mu_{d,n})=\mu_{d,n}$;\\
\> if $(d,n)$ is marked semi-critical, then $\II(\mu_{d,n})<2\mu_{d,n}$;\\
\> if $(d,n)$ is marked critical, then $\II(\mu_{d,n})=2\mu_{d,n}$.
\end{tabbing}
\item Let position $(d_0,n)$ be marked as either critical or semi-critical. Draw a line
north-east from this position, and do the same
from the position $(d_0,0)$ in the tableau. Then the marking
above the second line must be copied above the first one.

\item Suppose that $(d,0)$ is marked critical, $(d-k,k)$ is also critical, and
$(d-i,i)$ is off-critical for $i<k$. Assume that $(d,n)$ is semi-critical for
some $n$. Then $(d-k,n+k)$ is also semi-critical.

\end{itemize}

Finally, we say that a tableau is {\it recurrent} if
$$\sup \{d|\; (d,k)\text{ is critical for some }k>0\}=\infty;$$
we say that it is {\it periodic} if there exists $k>0$ such that the $k$-th column is entirely critical.
\end{defn}

\noindent
The relevance to the quadratic Yoccoz' puzzle should be evident from the above discussion:

\begin{defn}[{\bf The critical tableau of a Yoccoz' puzzle}] For $f_c$ as above, we let
$$\mu_{d,n}=\mod A_d(f_c^n(0)).$$
\end{defn}

We note:

\begin{prop}
The critical tableau of the Yoccoz' puzzle of $f_c$ is periodic if and only if $f_c$ is renormalizable.
\end{prop}

\noindent
The basis of the Yoccoz' result is given by the following theorem:

\begin{thm}
\label{tableau diverges}
Assume that $(\mu_{d,n})$ is a tableau, which is recurrent and not periodic. Assume further that
there exists $d$ such that $\mu_{d,0}>0$. Then
$$\sum_d \mu_{d,0}=\infty.$$

\end{thm}

\subsection{A puzzle partition for $R_a$.}
The puzzle pieces for $R_a$ which we construct are similar to those
just described but instead of external rays we use  bubble rays.
More specifically, choose a parameter $a$ in a parabubble wake
$W(t^+,t^-)$, and let the corresponding orbit portrait be
$$\OO(t^+,t^-)=\{\{\th_1,\ldots,\th_q \}\}.$$
Denote $\BB_i=\BB_{\th_i}$ the bubble ray with angle $\th_i$ starting with the bubble $A_\infty$,
and let $\al_a$ be the common landing point of these rays. Another repelling fixed point of $R_a$,
that in the intersection of $\bar A_0$ and $\bar A_\infty$ will be denoted $p_a$.

\begin{defn}

The {\it thin initial puzzle-pieces} of $R_a$ are the connected components of
$$\hat \C\setminus \left(\oli{(\cup_i \BB_i)}\cup\{\al_a\}\right).$$
Similarly, a {\it thick initial puzzle-piece} of $R_a$ corresponding
to a thin puzzle-piece $P$ is the
set
$$\bar P\cup \BB^1\cup \BB^2,$$
where $\BB^i$ are the two bubble rays which bound $P$.

Finally, an {\it initial puzzle-piece} of $R_a$ is a domain obtained
as follows. Let $\gamma_i$ be the axis of $\BB_i$ terminating at
$\al$ and $\infty$. Further, let
\[
\Phi: A_{\infty} \mapsto \hat{\C}  \sm \D
\]
be the B\"ottcher coordinate,
fix an arbitrary $r>1$, and let
$$D=\Phi^{-1}(\{|z|>r\})\text{ and }D'=R_a^{-1}(D_r)\cap A_0.$$
%
%
The initial puzzle-pieces are the connected components of
$$\hat\C\setminus\left( (\cup \gamma_i)\cup\{\al_a\}\cup\bar D\cup\bar D'  \right).$$
We denote the initial puzzle-pieces $P_0^1,\ldots,P_0^q$. The puzzle pieces of
{\it depth} $n$ are the $n$-th preimages of $P_0^i$, they will be denoted $P_n^j$.
\end{defn}

The basic properties being the same for all three kinds of puzzle-pieces, we will
only formulate the results for the last kind. We begin by noting:

\begin{lem}[{\bf Markov property}]
\label{markov}
For any two puzzle pieces $P_n^i$, $P_m^j$ one of the following two possibilities holds:
they are disjoint, or one is a subset of the other.
\end{lem}

\noindent
This allows us again to define for a point $z\in J(R_a)$ which is not a preimage of $\alpha_a$
$P_d(z)$ as the puzzle-piece of depth $d$ which contains $z$. Further, set
$$A_d(z)=P_d(z)\setminus \overline{P_{d+1}(z)};$$
we refer to this set as a complementary annulus, although it could  be degenerate.
We again label the annuli as critical, off-critical, and semi-critical depending on the
position of the critical point $-1$.
A critical annulus $A_{d+k}(-1)$ will be  called a {\em child} of the
critical annulus $A_d(-1)$ if
$$R_a^k:A_{d+k}(-1)\to A_d(-1)$$
is an unramified double covering.

We define $\TT_a$ to be a marked array
$$\TT_a=(\mod A_d(R_a^n(-1))),\;d,n\geq 0,$$
with the positions marked as critical, off-critical, or semi-critical if the
respective annuli are.
The following Proposition is verified in a straightforward way, completely similarly to the
quadratic case. We therefore omit the proof.

\begin{prop}
The marked array $\TT_a$ is a critical tableau.
\end{prop}

\noindent
However, it may happen that there is no non-degenerate annulus in the tableau $\TT_a$. We will
need to modify the construction of the annuli slightly to guarantee the existence of one.


\begin{figure}
\begin{center}
\includegraphics[scale=.9]{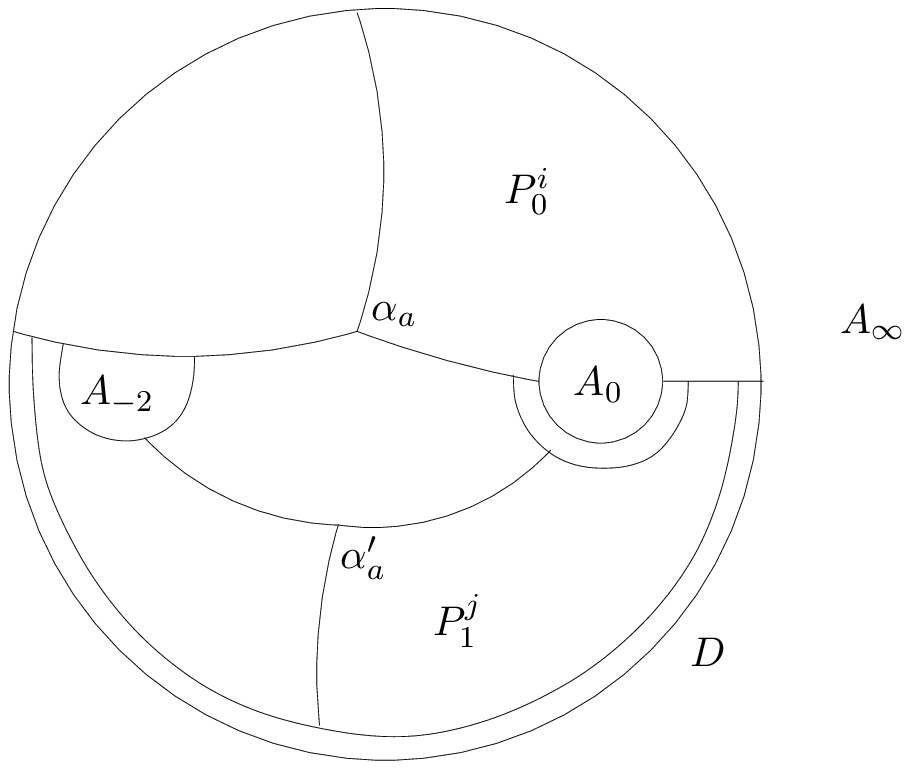}
\end{center}
\caption{A bubble puzzle of depth $1$. Note that the pieces $P_0^i$
  and $P_1^j$ touch at an arc connecting $A_0$ and
  $A_{\infty}$.\label{depth1-r} } 
\end{figure}


\begin{figure}

\begin{center}
\includegraphics[scale=1]{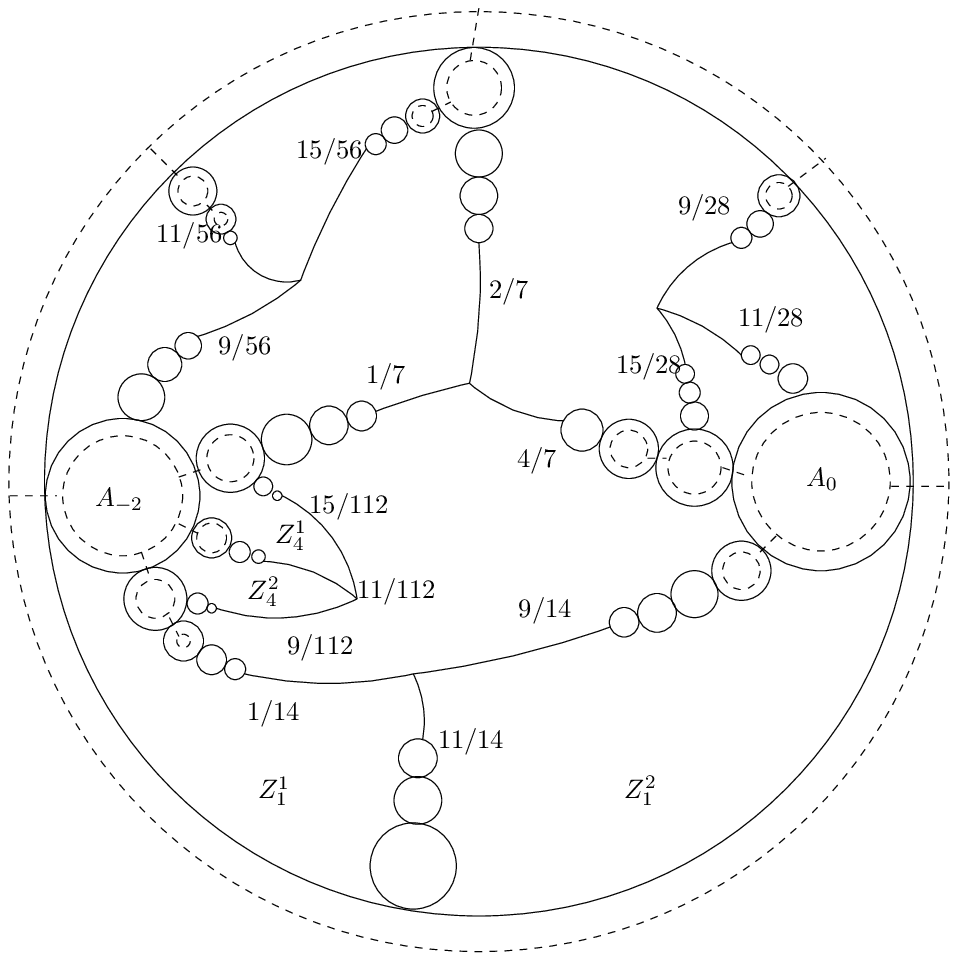}
\end{center}
\caption{A bubble-puzzle of  depth $1$ together with some preimages of
the pieces $Z_1^1$ and $Z_1^2$.
The broken lines
show the ``equipotential'' of depth $4$. Note
that $Z_4^1$ is degenerate in the sense that its boundary touches the
boundary of $P_0(-1)$, whereas $Z_4^2$ is not.\label{bpuzzle}} 
\end{figure}


\subsection{Non-degenerate annuli.}
The construction of a non-degenerate critical annulus for $R_a$ is somewhat
more delicate than that for a quadratic polynomial. We begin with the following:

\begin{Lem} \label{nondegI}
 We have $P_1(-1) \Subset \hat \C\setminus \bar D$.
\end{Lem}

\begin{proof}
There are $q \geq 3$ infinite bubble rays $\BB_k$, $k=1,\ldots
q$ landing at $\al$.
 First, let us argue that at least one bubble ray
$\BB_k$  contains $A_{-2}$ (the Fatou component of $R_a$ containing $-2$)
and another contains $A_0$. Suppose this is not the case. Then all,
but possibly one,
external angles $\th_k$ for $B_k$ will belong to $(1/6,1/3) \cup
(2/3,5/6)$. But then all, but possibly one, of the images of $\th_k$
under  doubling will belong to $(1/3,2/3)$, which
is disjoint from $(1/6,1/3) \cup (2/3,5/6)$. Since $q \geq 3$ this
gives a contradiction.

We want to show that the preimages
$\BB_k'$ of $\BB_k$ landing at the preimage of  $\al$ have
the same property, that is at least one bubble ray $\BB_k'$
contains $A_{-2}$ and another contains $A_0$. If this is not the case then
the images of all, but possibly one, $\BB_k'$ have angles in $(1/3,2/3)$,
which is impossible.

Hence the region $P_1(-1)$ is bounded by four bubble rays which all
emerge from $A_{-2}$ or $A_0$. It is easy to see that this region is
compactly contained in $A_{\infty}^c$, and the lemma follows.
\end{proof}



\comm{

Douady and Hubbard's definition of immediate renormalisation (see
e.g. \cite{Milnor5}) can be stated as a
combinatorial statement for $R_a$. If a quadratic polynomial $f_c$ is
DH renormalisable it is renormalisable by a ``thickening procedure''.
\begin{defn}
Assume that the initial puzzle for $R_a$ exists as above.
We say that $R_a$ has {\em DH renormalisable combinatorics} if and
only if the critical point $-1$ is mapped into $Y_0^1$ under
$R_a^{qn}$ for all $n \geq 1$. In other words, the critical point
never escapes $P_1(-1)$ under iterats of $R_a^q$.
\end{defn}

}

Now let us denote $Z^1_1,\ldots,Z^1_{q-1}$ the puzzle-pieces of level $1$ which
are not adjacent to $\al_a$, but to its other preimage $\al_a'$.
It is easy to see that if $\TT_a$ is not a periodic tableau, then some iterate
of the critical point $-1$ under $R^q_a$ will escape to one of the
pieces $Z_1^j$.
%
The first time this happens, say after the $n$-th iterate, we can pull
back the {\em degenerate} annulus $P_0(-1) \sm P_1^j$ under
$R_a^{qn}$.
See Figure \ref{bpuzzle} for an illustration.
This will give a degenerate critical annulus
$A_m(-1)$. However, by \lemref{nondegI}, the only place where the boundaries of  $P_m(-1) $
and $P_{m+1}(-1)$ touch is a preimage of the segment $l$ of two internal rays containing
$\oli{A_{\infty}} \cap \oli{A_0}$  which connects $D$ and $D'$. The invariance of $A_\infty\cup A_0$ implies:

\begin{lem}
\label{pinch}
The pinching of any child of $A_m(-1)$ is disjoint from
  $\overline{P_m(-1)}\cap \overline{P_{m+1}(-1)}$.
\end{lem}

This means in particular the following:
\begin{cor}  \label{thick}
Let $A_m(-1)$ be as above.
 Let $A_{m_j}(-1)$ be any child of
$A_m(-1)$. Then the critical puzzle pieces $P_{m_j}(-1)$ satisfy
\[
P_{m_{j+1}}(-1) \Subset P_{m_j}(-1)  \text{ and }
P_{m_{j+1}+1}(-1) \Subset P_{m_j+1}(-1).
\]
\end{cor}

\comm{
\noindent
Let us now {\it thicken} the arc $l$ to a strip $S$ which is invariant under the
branch of $R_a^{-1}$ fixing $p_a$. Let us now replace $A_m(-1)$ with a non-degenerate
annulus $\tilde A_m(-1)$ obtained by thickening $A_m(-1)$ with the preimage of $S$.

\begin{defn}
\label{thick tableau}
The {\it thickened} tableau $\tilde\TT_a$ is obtained from $\TT_a$ by preserving the
same marking, and replacing the moduli of the images and pre-images of $A_m(-1)$ with those
of $\tilde A_m(-1)$.
\end{defn}
}

\noindent
We first handle the non-recurrent case:
\begin{Lem} \label{non-recurrent}
If there is some $N$ so that $P_N(-1)$ is disjoint from the orbit $z_0
\mapsto z_1 \mapsto \ldots$, then $\cap_n P_n(z_0) = \{z_0 \}$.
\end{Lem}
\begin{proof}
The proof goes as in \cite{Milnor5}. We first thicken the puzzle-pieces of level $N-1$
to domains
$U_i \supset{P}^i_{N-1}$,
numbered so that $U_0\supset P_{N-1}(-a)$ and with the following property: for each $i>0$ there are 
 two univalent branches $g_1^i$ and $g_2^i$ of $R_a^{-1}$ defined on $U_i$, each of which carries it
into a proper subset of some $U_j$.
This is easily done, we leave the details to the reader.
We next equip every $U_i$ with the Poincar\'e distance $\rho_i(x,y)$.
It follows that for each puzzle piece $P_{N-1}^i$, $i>0$,
the branch $g_k^i$ shrinks the Poincar{\'e} distance by some definite
factor $\lambda<1$.
Since the orbit $z_0,z_1,z_2,\ldots$ avoids the
critical puzzle piece we get that
\[
\diam(P_{N+h}(z_0)) \leq \de \la^h,
\]
and the statement of the lemma follows.

\end{proof}

We next attack the harder recurrent case:
\begin{thm}
\footnote{We thank Carsten Petersen for pointing out that the proof of \thmref{pieces shrink} given in the published version
may fail in some cases. We supply the corrected proof below.}

\label{pieces shrink}
Assume that the critical tableau $\TT_a$ is recurrent and not periodic.
Then
$$\bigcap P_d(-1)=\{-1\}.$$

\end{thm}

Assume $A$ is a degenerate critical annulus. We may assume that $A$ is excellent. 
Hence every child is excellent as well.
This forms a tree of descendants $A_{i,j}$ starting from $A=A_{0,1}$ so that, for fixed $i > 0$, 
$A_{i,j}$ are the descendants of generation $i$. Generation $i$ means that $f^i(A_{i,j})=A_0$ and 
that $f^k: A_{i,j} \raw A_0$ is a $2^i$ degree unbranched covering. Moreover, since every
 $A_{i,j}$ is excellent there are at least $2^i$ annuli of generation $i$.

All $A_{i,j}$ form a nest around the critical point. We can relabel them so that 
$A=A_0$ surrounds $A_1$ which in turn surrounds $A_2$ and so on. In this way 
we get a nested sequence of annuli.

The {\em complementary annulus} $\al_j$ is defined to be the  
annulus between the outer boundary of  $A_j$ and the outer boundary of $A_{j-1}$. 
Note that by Corollary \ref{thick} any complementary annulus is non-degenerate.


\thmref{pieces shrink} will follow from:

\begin{Lem} \label{infinitemod}
The sum of the moduli of all complementary annuli is infinite.
\end{Lem}

An annulus is a difference between two puzzles pieces, an outer puzzle piece and an inner puzzle piece, 
the inner being contained in the outer. We say that an annulus $A$ {\em surrounds} a set $E$ if 
the inner puzzle piece of $A$ contains $E$.

Take some complementary $\al$ which lies between the two degenerate annuli $P=A_l$ and $Q=A_{l+1}$, where 
$P$ surrounds $Q$. Note that we assume that no annulus $A_j$ lies strictly between $P$ and $Q$. 
Now $Q$ has a child, say $Q_1$, so that $Q_1$ maps onto $Q$ as a $2$ degree unbranched 
covering. We want to pull back $P$ along the same branch (if possible) as $Q$ back to some $P_j$ surrounding $Q_1$.

In the first steps $\al$ (between $P$ and $Q$) is pulled back as a one-to-one map 
until some preimage $P_1$ of $P$ under $f^{k}$ surrounds the critical point.
This means by definition that this preimage $P_1$ is a child to $P$.
 If, moreover, $Q_1$, being the preimage of $Q$ under $f^{k}$ surrounded by $P_1$, also surrounds the critical point, then
 we stop and have found $P_1$ surrounding $Q_1$ both being children of $P$ and $Q$ respectively.
 Since we assumed that no degenerate annulus $A_j$ is between $P$ and $Q$, it follows that 
there cannot be any such $A_i$ between $P_1$ and $Q_1$ either.

The second, and more likely, case is that, whereas $P_1$ surrounds the critical point, $Q_1$ does not surround the critical point. 
Hence we are in a semi-critical situation, so the pullback $f^{-k}(\al)$ is not an annulus. 
However, if we consider the annulus $\be_1$ between $P_1$ and $Q_1$, this annulus has modulus at
 least $1/2$ of the modulus of $\al$, by a standard inspection from semi-critical annuli. 
Continuing pulling back $\be_1$, we again sooner or less reach the same situation:
Some pullback $P_2$ of $P_1$ under $f^{k_1}$ surrounds the critical point. If again the preimage $Q_2$ (being a preimage of $Q_1$ under $f^{k_1}$) surrounded by $P_2$ also surrounds the critical point we are done and have found two descendants $P_2$ and $Q_2$ to $P$ and $Q$ respectively. However, note that, whereas $Q_2$ is a child to $Q$, we have that $P_2$ is a child of $P_1$ and $P_1$ is a child of $Q$. 
$Q_1$ is not a child of $Q$ since $Q$ was assumed to be disjoint from the critical point.

Continuing in this way we find two descendants $P_m$ and $Q_m$ such that
$$f^{k+k_1+\ldots+k_{m-1}}: P_m \raw P$$
as a $2^m$ degree unbranched covering and
$$f^{k+k_1+\ldots+k_{m-1}}: Q_m \raw Q$$
as a $2$ degree unbranched covering.

Hence, $Q_m$ is a child to $Q$, whereas every $P_{j+1}$ is a child to
$P_j$, $j=0, \ldots, m-1$. In this case we call the annulus between
$P_m$ and $Q_m$ an {\em offspring} of $\al$. 
Hence, every offspring has modulus at least $2^{-m}$ times the modulus
of its {\em ancestor} $\al$, where $m$ is as defined above. 
Again, there cannot be any degenerate annulus $A_j$ between $P_m$ and $Q_m$. Otherwise, we could map 
this annulus forward: $f^{k+k_1+\ldots+k_{m-1}}(A_j)$ would be a degenerate annulus between $P$ and $Q$.

Conversely, let $P_m$ and $Q_m$ be given degenerate annuli surrounding the complementary annulus $\al_1$ and assume that there is no other degenerate annulus between $P_m$ and $Q_m$. If $Q_m$ has generation more than $1$ then the parent $Q$ would have generation at most $1$. On the other hand, the parent $P$ to $P_1$, which in turn is parent to $P_2$ and so on down to $P_m$, might have negative generation, meaning that $P$ is actually a parent to $A_0$. In this case, $A_0$ would lie between $P$ and $Q$. But in this case there has to be some preimage of $A_0$ laying between $P_m$ and $Q_m$. This contradict the fact that there is no degenerate annulus between $P_m$ and $Q_m$.

We conclude from the above discussion:
\begin{Lem}
Every complementary annulus $\al$ between two degenerate annuli $P$ and $Q$, where the generation of $Q$ is larger than $1$, 
has some unique ancestor $\be$.
\end{Lem}

\begin{Def}
Given a complementary annulus $\al$ surrounded by the outer degenerate annulus $A_{m,*}$ and the inner annulus $A_{n,*}$, we say that the {\em outer generation} to $\al$ is equal to $m$ and the {\em inner generation} to $\al$ is $m$. We write $\al=\al_{n,*}^m$, where $*$ means an index, since there might be many $\al$ with the same $m$ and $n$.
\end{Def}

We have proved the following.

\begin{Lem} \label{onestep}
For every complementary annulus $\al=\al_{n,*}^m$ with $n > 1$ and with ancestor $\al_{n-1,*}^{m_1}$ we have
\[
\mod(\al_{n,*}^m) \geq 2^{m_1-m} \mod(\al_{n-1,*}^{m_1}).
\]
\end{Lem}

\begin{Cor} \label{manysteps}
For every complementray annulus $\al_{n,*}^m$, $n > 1$, there is some ``grand'' ancestor $\al_{1,*}^{m_{n-1}}$ such that
\[
\mod (\al_{n,*}^m) \geq 2^{m_{n-1}-m} \mod (\al_{1,*}^{m_{n-1}}).
\]
\end{Cor}

Since the number of degenerate annuli of generation $m$ is at least $2^m$ we have that the number of complementary annuli of outer generation $m$ is at least $2^m$.
Moreover, trivially, we have $\mod(\al_{1,*}^m) \geq M_0$ for all $m$,
for some $M_0 > 0$. 

By Corollary \ref{manysteps} the sum of the moduli of all $\al_{n,*}^m$ for fixed $m$ is at least
\[
\sum_{n,*} \mod (\al_{n,*}^m) \geq 2^m 2^{-m} \mod (\al_{1,*}^{m_{n-1}}) \geq M_0.
\]
Hence
\[
\sum_{m,n,*} \mod( \al_{n,*}^m ) = \infty,
\]
and Lemma \ref{infinitemod} follows.

\subsection{Combinatorics of the puzzle.}
We make some definitions first. Let $a_1$, $a_2$ be two parameters in the
same wake $W$. We say that $R_{a_1}$ and $R_{a_2}$ have the {\it same
combinatorics of the puzzle up to depth} $d$ if there exists an
orientation preserving homeomorphism $\phi:\hat \C\to\hat\C$ such that the
following holds:
\begin{itemize}
\item $\phi$ homeomorphically maps distinct puzzle pieces $P_k^i$ of depth $k\leq d$
of $R_{a_1}$ to distinct puzzle-pieces $Q_k^j$ of depth $k$ of $R_{a_2}$;
\item for all $k\leq d$ we have $\phi:P_k(-1)\mapsto Q_k(-1)$;
\item finally, $\phi$ {\it respects the dynamics}, that is,
$$P_k^i=R_{a_1}(P_k^j)\text{ if and only if }\phi(P_k^i)=R_{a_2}(\phi(P_k^j)).$$
\end{itemize}

Similarly, we will say that a quadratic polynomial $f_c$ and $R_a$ have
the same combinatorics of the puzzle up to depth $d$, if there exists an orientation-preserving
continuous surjection $\phi$ which maps puzzle-pieces of $f_c$ to those of $R_a$ up to depth $d$,
sending critical pieces to critical ones, and respecting the dynamics.


\begin{prop}
\label{parapiece}
Let $f_c$ be a quadratic polynomial without non-repelling fixed points. For every
$d$ there exists a parameter $a$ such that $R_a$ and $f_c$ have the same combinatorics
of the puzzle down to depth $d$. Moreover, consider the puzzle-piece $P_d(c)$ of $f_c$,
and let $\be_1,\ldots,\be_k$ be the angles of external rays which bound it. Then
bubble rays with the same angles bound the puzzle piece $Q_d(-a)$ of $R_a$.

Finally, there exists an open set $\Delta_d$ in the $a-$plane, with $\Delta_d\subset \Delta_{d-1}$,
and $\Delta_0=W$ such that $R_b$ has the same combinatorics of the puzzle to depth $d$
and $-b$ is contained in the particular puzzle piece of level $d$ if and only if
$b\in\Delta_d$.

\end{prop}
\begin{proof}
The Proposition follows by a straightforward induction on the depth
$d$. The base of induction, with $d=0$ is given by \lemref{wakesII}.
Assuming the statement is true at depth $d-1$, consider the pullback
of the puzzle of level $1$ inside the critical value piece
$P_{d-1}(-a)$. By assumption, this picture has the same
combinatorial structure as the similar one for $f_c$. By
\lemref{orbit persist}, as the parameter $a$ moves through
$\Delta_{d-1}$, the critical value sweeps out $P_{d-1}(-a)$. We can
hence select a parameter $a$ to match the combinatorics of the
puzzle of $f_c$ down to level $d$. The parameter plane statement
follows from similarly obvious consideration and is left to the
reader.

\end{proof}

\begin{defn}
We call a set $\Delta_d$ as above {\it a parameter puzzle piece}.
\end{defn}

\section{Existence of a Mating}

\comm{

Similarily to $f_c$, we get at every level $n$ of the Yoccoz
bubble-puzzle for $R_a$ we get a compact parameter set $\De_n$,
for which the critical truncated symbol sequence $\si(-1)$ of level
$n$ is the same critical truncated symbol sequence $\iota(0)$
for $f_c$. Note that $\De_{n+1} \subset \De_n$. Taking $a \in \cap_n
\De_n$, we get a rational map $R_a$ with the same Yoccoz
puzzle at any depth. In particular, the Yoccoz bubble-puzzle for $R_a$
is not periodic, since $f_c$ is non-renormalisable, by Lemma \ref{period-ren}.
We conclude that the critical puzzle pieces for $R_a$
also shrinks to a point, namely the critical point $-1$, by
Proposition \ref{mprop}.  Let $P_n^*(x)$ be the (finite number of) puzzle pieces which
intersect $x$. We immediately get the following.

\begin{Prop} \label{localc}
There is a parameter $a \in \mat$ such that
the critical tableau for $R_a$ is the same as the critical tableau for
$f_c$. For such $a$, all puzzle pieces $P_n^*(x)$, for $x \in J(R_a)$
shrink to a single point, i.e. $\cap_n P_n^*(x) = \{x \}$.
\end{Prop}

Since backward iterates of $-1$ is dense on the Julia set, and the
fact that bubble rays cannot cross each other, we conclude the following.

\begin{Cor} \label{corre}
Every infinite bubble ray lands at a single point.
The Julia set $J(R_a)$ of $R_a$ is locally connected.
\end{Cor}

}
Fix a Yoccoz' polynomial $f_c$ which is not critically finite, non-renormalizable,
and such that $c$ does not belong to the $1/2$-limb of the Mandelbrot set.
By \propref{parapiece}, there exists a parameter value $a$
such that $R_a$ has the same combinatorics of the puzzle as $f_c$
for all $d\in\N$.

\begin{lem}
\label{pieces}
Consider any $z\in J(R_a)$ which is not a preimage of $\al_a$ or $p_a$.
Then the nested sequence of puzzle pieces $P_d(z)$ shrinks to $z$:
$$\bigcap P_d(z)=\{z\}.$$

\end{lem}
\begin{proof}
Assume first that there exists some $N>0$
such that the orbit of $z$ is disjoint from $P_N(-1)$. In this case, the claim is implied
by \lemref{non-recurrent}.

In the opposite case, for each $n\in\N$ consider the first instance $i$ such that
$R_a^i(z)\in P_{n+1}(-1)$. Then the complementary annulus $\alpha_{n+i}(z)$ is a conformal copy of $\alpha_n(-1)$.
By construction, all these annuli around $z$ are disjoint, and hence by \lemref{infinitemod},

$$\sum \mod \alpha_n(z)=\infty.$$

\noindent
By \lemref{moduli diverge}, we have the claim.

\end{proof}

\begin{lem}
\label{rays land}
Every bubble ray for $R_a$ lands.

\end{lem}
\begin{proof}
This is obviously true for the preimages of the rays landing at the fixed point $\al$.
Let $z$ be an accumulation point of any other ray $\BB=\cup_0^\infty F_i$. There is an infinite sequence of
nested puzzle pieces $P_d(z)$ containing $z$, and by the previous Lemma,
$$\bigcap P_d(z)=\{z\}.$$
Now by \lemref{Ra} the bubbles $F_i$ do not cross the boundaries of $P_d(z)$, and hence $$F_i\to z.$$
\end{proof}


\subsection{Construction of semiconjugacies}
Consider the conjugacy
$$\phi:\cir{K}_{\bas} \cup_n f_{\bas}^{-n}(\al) \mapsto \cup_n R_a^{-n}(A_\infty\cup \{p\})$$ defined in Proposition
\ref{marking-mated}. By \lemref{rays land} and \lemref{bubble lc} it extends by continuity to
a semi-conjugacy $K_\bas\to \overline{\cup R_a^{-n}(A_\infty)}=\hat\C$:
\begin{equation} \label{phi1}
\phi_1 \circ f_{\bas} (z) = R_a \circ \phi_1 (z).
\end{equation}

Let $z\in J_c$ and not a preimage of $\al$,
and let $P_d(z)$ be the sequence of Yoccoz' puzzle-pieces of depth $d$ containing $z$.
Let $Q_d(z)$ be the corresponding pieces in the puzzle of $R_a$ and define
$$\phi_2(z)=\cap Q_d(z).$$
By construction, $\phi_2$ extends continuously to $\cup_n f_c^{-n}(\{\al\})$ and for the extended map
$$\phi_2\circ f_c=R_a\circ \phi_2.$$

\noindent
Let $\sim_r$ denote the ray equivalence relation generated by the quadratics $f_\bas$ and $f_c$.
We proceed to demonstrate:

\begin{thm}
\label{ray equivalence}
We have
$$\phi_i(z)=\phi_j(w)$$
if and only if they are in the same ray equivalence class,
$$z\sim_r w.$$
\end{thm}

\noindent
We begin with the following definition.

\begin{defn}
For $q>1$, let
$$\th_1\mapsto\th_2\mapsto\cdots\mapsto\th_q\mapsto\th_1$$
be a period $q$ orbit of the doubling map. The angles $\th_i$ partition
the circle into arcs $A_i$, $i=1,\ldots,q$, which we enumerate in the counter-clockwise
order starting from the arc containing $0$. For $\th\in \T$ which does not eventually fall
into the orbit under doubling, we denote
$\sigma_{\th_1,\ldots,\th_q}(\th)$ the {\it itinerary} of $\th$ with respect to the partition
$A_i$, viewed as an infinite string in $\{1,\ldots,q\}^\infty.$ In the case when $\th$ is a
preimage of one of the $\th_i$ the itinerary $\sigma_{\th_1,\ldots,\th_q}(\th)$ will be a finite
string of digits between $1$ and $q$ -- to avoid ambiguity, the last $A_i$ will be chosen to the
right of $\th_i$.
\end{defn}

\comm{

\noindent Given the expansiveness of $d:z\mapsto 2z \mod \Z$, we
have:

\begin{prop}
\label{same itinieraries}
For $\xi,\eta\in\T$,
$$\sigma_{\th_1,\ldots,\th_q}(\xi)=\sigma_{\th_1,\ldots,\th_q}(\eta)\text{ if and only if }\xi=\eta.$$

\end{prop}
}


In a very similar way, let us define a symbol sequence $\si(z) \in
\{1,\ldots,q\}^{\infty}$ with respect to the initial Yoccoz puzzle
for $f_c$ or the initial Yoccoz bubble-puzzle for $R_a$ as follows.
Enumerate the initial puzzle-pieces of $f_c$ as  $P_0^k$, $k=1,\ldots,q$ in counter-clockwise
order around $\alpha$, starting with $P_0^1\ni 0$. Set $Q_0^k$ to be the puzzle piece
of $R_a$, which corresponds to $P_0^k$.
Put

\begin{equation}
\si(z) = \begin{cases} k & \text{if $f_c^j(z) \in P_0^k$, for $z \in
J(f_c) \sm \cup_n f_c^{-n}(\al)$},
  \\
k&   \text{if $R_a^j(z)\in Q_0^k$, for $z \in J(R_a) \sm \cup_n
R_a^{-n}(\al_a\cup p_a)$}.
\end{cases}
\nonumber
\end{equation}


Since $\phi_1$ is a semi-conjugacy the following lemma is immediate.
\begin{Lem} \label{sigma}
Assume that $z \in K_{\bas}$ is uni-accessible and let
$\phi_1(z)=\zeta$. Let $-\be$ be the angle of the external ray
landing at $z$. If $z$ is not a preimage of the $\alpha$-fixed point, then
\[
\si(\zeta)=\si_{-\th_1,\ldots,-\th_q}(-\be).
\]
\end{Lem}

\noindent Recall now, that a point in the Julia set $J_\bas$ is
bi-accessible if and only if it is a preimage of $\alpha_\bas$. The
latter is the landing point of two external rays, $R_{1/3}$, and
$R_{2/3}$, forming a period $2$ cycle. \noindent Let $d$ be the
function $d:z\mapsto 2z \mod \Z$.

\begin{lem}
\label{ray eq}
Let $R_\th$ be a ray landing at a bi-accessible point $x\in J_\bas$. Then the landing point of
$R_{-\th}$ in $J_c$ is uni-accessible.
\end{lem}

\begin{proof}
The angle $-\th$ has a finite orbit under the doubling, and hence the orbit of the landing point $y$ of the ray
$R_{-\th}$ is also finite. By assumption, $f_c$ is not critically finite, and hence the orbit of $y$
does not include $0$. Denote $n$ the first iterate for which $d^n(-\th)\in\{1/3,2/3\}$,
and $z=f^n(y)$. Since $f^n$ is a local homeomorphism on a neighborhood of $y$, the number $m$
of accesses is the same for $y$ and $z$. Assume that $m>1$.

Note first that $z$ cannot be a fixed point, as otherwise the ray portrait $\{\{1/3,2/3\}\}$ is
realized for $f_c$, and $c$ is in the $1/2$-limb. Hence $z$ has period $2$.
By the properties of periodic external rays all rays landing at $z$ have the same period, $2$,
and same for $f(z)$. Hence, there are $m\times 2\geq 4$ angles in $\T$ whose period under the doubling is
equal to $2$. By inspection, $1/3$ and $2/3$ are the only angles with this property, and we have
arrived at a contradiction.
\end{proof}

\noindent
By assumption, there exists $q>2$ such that there is a cycle of rays $R_{\th_1},\ldots,R_{\th_q}$
landing at the dividing fixed point $\al$ of $f_c$. By construction, a cycle of bubble rays
$\BB_{\th_1},\ldots,\BB_{\th_q}$ with the same angles lands at the fixed point $\al_a$.

\begin{lem}
\label{phi1 is ok}
We have
$$\phi_1(z)=\phi_1(w)\text{ if and only if }z\sim_r w.$$
\end{lem}
\begin{proof}

\comm{ By \lemref{Ra}, both $z$ and $w$ are landing points of
infinite bubble rays $\BB_1,\BB_2\subset K_\bas$. Denote
$\be_i=\angl(\BB_i)$, and let $\BB'_i=\phi_1(\BB_i)$. Note that by
definition, the angles of the external rays of $f_\bas$ landing at
$z$ and $w$ are $-\be_1$ and $-\be_2$ respectively.

Consider the external rays $R_{\be_i}$ of $f_c$. Since the
combinatorics of the puzzles of $f_c$ and $R_a$ is the same for
every depth, these two rays have a common landing point. The
statement of the lemma follows from \lemref{ray eq}. }


By \lemref{Ra}, only uni-accessible points can be identified. From
Lemma \ref{ray eq} the lemma now follows if at least one of $z$ and
$w$ is bi-accessible. Hence we can assume that both $z$ and $w$ are
either landing points of infinite bubble rays $\BB_1,\BB_2\subset
K_\bas$, or that one of $z$ and $w$ or both lies on a uni-accessible
point on the boundary of a bubble. Denote $-\be_1,-\be_2$ the angles
of the external rays landing at $z$ and $w$ respectively.
(In the case when $z$ and $w$ are landing points of infinite bubble
rays $\BB_i$, note by definition, that the angles of these bubbles
rays are $\be_1$ and $\be_2$ respectively.)

By Lemma \ref{sigma}, $\phi_1(z)=\phi_1(w)$ if and only if
\begin{equation} \label{sieq}
\si_{-\th_1,\ldots,-\th_q}(-\be_1)=\si_{-\th_1,\ldots,-\th_q}(-\be_2).
\end{equation}

Now, consider the external rays $R_{\be_i}$ of $f_c$. Since the
combinatorics of the puzzles of $f_c$ and $R_a$ is the same for
every depth, these two rays have a common landing point if and only
if (\ref{sieq}) holds. The statement of the lemma now follows from
\lemref{ray eq}.
\end{proof}

\begin{lem}
\label{phi2 is ok}
We have
$$\phi_2(z)=\phi_2(w)\text{ if and only if }z\sim_r w.$$
\end{lem}
\begin{proof}
Note that  by \lemref{ray eq}, if $z\neq w$, then
$z\sim_r w$  if and only if both of these points are uni-accessible, and denoting
$\be_1$, $\be_2$ their external angles, we have $d^n(\be_1)=1/3$, $d^n(\be_2)=2/3$ for some $n$.

On the other hand, if $\zeta=\phi_2(z)=\phi_2(w)$, then $\zeta\in R_a^{-n}(p_a)$ for some $n$.

It is thus enough to show, that $\phi_2(z)=\phi_2(w)=p_a$ if and only if $z$, $w$ are the landing points
of the external rays $R_{1/3}$, $R_{2/3}$ respectively. By construction, at most two points in $J_c$
are mapped to $p_a$ by $\phi_2$, so we only need to prove the second implication.

The landing points $z$, $w$ of rays $R_{1/3}$, $R_{2/3}$ form a cycle of period $2$, hence, the period
of the cycle $\zeta_1=\phi_2(z)$, $\zeta_2=\phi_2(w)$ is at most $2$. By \lemref{Ra}, these points do not lie in the
boundary of any of the bubbles. Assume that $\zeta_1\neq p_a\neq \zeta_2$. Then there exists a bubble ray of
angle $\th$ landing at $\zeta_1$. Since the combinatorics of the puzzle is the same for $R_a$ and $f_c$,
$$\sigma_{\th_1,\ldots,\th_q}(\th)=\sigma_{\th_1,\ldots,\th_q}(1/3).$$
This bubble ray then lands at a point in $J_\bas$ with the external
angle $2/3$, which is a contradiction.

\end{proof}

We finish the proof of \thmref{ray equivalence} with the following:

\begin{lem}
\label{both phi12 ok}
We have $\phi_1(z)=\phi_2(w)$ if and only if $z\sim_r w$.
\end{lem}
\begin{proof}

\comm{

Let us prove this statement first in the case when $z$ is an
endpoint of an infinite bubble ray $\BB\subset K_\bas$. Let
$\BB'=\phi_1(\BB)$ land at $\zeta\in J(R_a)$. Let $\be$ be the
external angle of the point $z$; by definition,
$$\angl(\BB)=\angl(\BB')=-\be.$$

If $\zeta=\phi_2(w)$, then consider the external ray $R_\gamma$
landing at $w$, which passes through the same puzzle-pieces, as the
bubble-ray $\BB'$. We have
$$\sigma_{\th_1,\ldots,\th_q}(\gamma)=\sigma_{\th_1,\ldots,\th_q}(\angl(\BB')),$$
and hence $\gamma=\angl(\BB')$, and $z\sim_r w$.

The ``if'' direction, when $z$ is a landing point of a bubble ray is
evident by construction. The general case follows by continuity. }


If $z \in K_{\bas}$ is uni-accessible then let $-\be$ be the angle
of the external ray landing at $z$ and put $\zeta = \phi_1(z)$. By
Lemma \ref{sigma},
\[
\si_{-\th_1,\ldots,-\th_q}(-\be)=\si(\zeta).
\]

If $\zeta=\phi_2(w)$, then $w$ lies in the same puzzle-pieces as the
point $\zeta$, by definition. An external ray $R_\gamma$ (there can
be more than one) which lands at $w$ must by Lemma \ref{sigma}
satisfy
$$\sigma_{\th_1,\ldots,\th_q}(\gamma)=\si(\zeta).$$
Obviously, one solution is $\gamma=-\be$, and therefore $z\sim_r w$.
Conversely, if $z\sim_r w$, then $\phi_1(z)=\phi_2(w)$ by
construction.

If $z \in K_{\bas}$ is bi-accessible then the lemma follows from
Lemma \ref{ray eq}.

\end{proof}

\noindent
We conclude:

\medskip
\noindent
{\bf Main Theorem, the existence part.} {\it Suppose $c$ is a non-renormalizable parameter value outside the $1/2$-limb of
$\MM$. Then the quadratic polynomials $f_c$ and $f_\bas$ are conformally mateable.}

\section{Uniqueness of mating}
 \label{parapuzzle}
To transfer the results of shrinking puzzle pieces in the dynamical
plane to the parameter plane, we use a variation of the
approach of Yoccoz (see \cite{Hubbard}). Our arguments follow the
presentation of \cite{PR}.

Let us recall the following definition.
\begin{Def}
Let $X$ be a connected complex mainfold.
A {\it holomorphic motion} over a set $E \subset \C$ is a function
\[
\varphi: X \times E \raw \hat{\C},
\]
where $\varphi (\la,z)$ is holomorphic in the variable $\la \in X$ and
injective in $z \in E$.
\end{Def}

\noindent
We make use of a stronger version of the $\la$-lemma of Man{\~e}-Sud-Sullivan \cite{MSS},
due to Slodkowski \cite{Slo}.

\begin{Lalemma}\label{McMullenbook}
A holomorphic motion over a set $E$ has a unique extension to a
holomorphic motion over $\oli{E}$. The extended motion gives a
contiuous map $\varphi: X \times \oli{E} \raw \hat{\C}$. For each
$\la \in X$, the map $\varphi_{\la}: \oli{E} \raw \hat{\C}$ extends to
a quasiconformal map of $\hat{\C}$ to itself.
\end{Lalemma}

\noindent
Let us fix a parameter $c$ satisfying the conditions of the Main Theorem. Let $\Delta_n$ be the nested
sequence of parameter puzzle-pieces of \propref{parapiece} in the $a$-plane. Our aim is to show:

\begin{thm}
\label{parapieces shrink}
We have $$\diam (\Delta_n)\to 0.$$
\end{thm}

\noindent Let us fix a parameter $a_0\in \cap \Delta_n$. Let $P$ be
a parabubble intersecting some $\Delta_n$. Denote $B_a$ the bubble
containing the critical value $-a$ for $R_a$ with $a\in P$. Let
$k\in\N$ be the smallest such that for any $a\in P$, $R^k(-a) \in
A_{\infty}^a$. Let
$$\Phi_a:A_\infty^a\to\hat\C\setminus\D$$
denote the normalized B{\"o}ttcher coordinate at infinity. By \lemref{bubble lc},
it extends homeomorphically to the boundary.
We then obtain a homeomorphism $P\mapsto B_{a_0}$ by the formula.
\[
F: a \mapsto R_{a_0}^{-k} \circ \Phi_{a_0}^{-1} \circ \Phi_a \circ
R_a^k(-a).
\]

\noindent
Pasting these homeomorphisms together, we obtain

\begin{Lem} \label{hbound}
There is a homeomorphism from the closure of the boundary of
the parameter puzzle piece of depth $n$ into the closure of the boundary of the
puzzle of depth $n$ for $R_{a_0}$.
\end{Lem}

We now construct a holomorphic motion on the boundary of the
puzzle at an initial level.

\begin{Lem} \label{initial}
There is a holomorphic motion $h_n: \De_n \times I_{n+1}^{a_0} \raw
\hat{\C}$, where $I_{n+1}^{a_0}$ is the closure of the boundary of
the puzzle of depth $n+1$. We have $h_n^a(I_{n+1}^{a_0}) =
I_{n+1}^a$. Moreover, $R_a \circ h_n^a (z) = h_n^{a_0} \circ
R_{a_0}(z)$, for any $z \in I_{n+1}^{a_0}$.
\end{Lem}
\begin{proof}
Indeed, as $a$ varies throughout $\De_n$, the critical value does
not hit the bubble rays corresponding to the puzzle of depth $n$
according to Lemma \ref{above}. We get from Lemma \ref{boundary}
that $A_{\infty}^a$ moves holomorphically on $\De_n$. So do the
preimages of $A_{\infty}^a$ as long as we do not hit the critical
value. It follows that every bubble $B$ in the boundary of the
puzzle of depth $n$ moves holomorphically according to the formula
\begin{equation}
\eta_a(z)=R_{a}^{-k} \circ \Phi_{a}^{-1} \circ \Phi_{a_0} \circ
R_{a_0}^k(z), \label{formula}
\end{equation}
where $k$ is smallest integer such that $R_a^k(z) \in A_{\infty}$,
for $z \in B$.

Since the critical value does not intersect the puzzle of depth $n$,
we can pull back this puzzle once so that the puzzle of depth $n+1$
moves holomorphically as well.

The $\la$-Lemma now extends the motion to its closure. It follows
from (\ref{formula}) that $h_n^a(I_{n+1}^a)=I_{n+1}^{a_0}$ and that
the diagram
\[
\begin{CD}
I_{n+1}^{a_0} @> {h_n^a}>> I_{n+1}^a          \\
@V {R_{a_0}} VV        @VV  {R_a}  V          \\
I_n^{a_0} @> {h_{n-1}^a}>> I_{n}^a
\end{CD}
\]
is commutative.
\end{proof}

\begin{defn}
Let $D_{n+1}^a$ be the puzzle piece bounded by $h_n^a(\partial
P_{n+1}^{a_0})$, where $P_{n+1}^{a_0}$ is the puzzle piece $P_{n+1}$
surrounding the critical value $-a_0$ at depth $n+1$.
\end{defn}

\noindent
We have the following:
\begin{Lem} \label{diff}
The parameter $a \in \De_m \sm \De_{m+1}$ if and only if the critical
value $-a \in D_m^a \sm D_{m+1}^a$.
\end{Lem}
\begin{proof}
Take a non self-intersecting path $a_t$ from $a_0$ to the boundary
of $\De_m$ ,$t \in [0,1]$, crossing the boundary of $\De_{m+1}$
exactly once. Then the critical value $-a_t$ has to cross the
boundary of $h_m^{a_t}(\partial P_{m+1}^{a_0})$, since we always
have $D_m^a \supset D_{m+1}^a$. Assume this happns at $t=t_0$. Then
for $t > t_0$ we get that $-a_t \notin D_{m+1}^{a_t}$, since we are
outside $\De_{m+1}$. Similarily, $-a_t \in D_{m+1}^{a_0}$ for $t <
t_0$.
\end{proof}


\begin{proof}[Proof of \thmref{parapieces shrink}]

Let us first handle the harder case, when the critical tableau of $f_c$ is recurrent.

Extend the holomorphic motion on $\De_{m_0}$ at depth $m_0$ by the
$\la$-Lemma, so that we get a holomorphic motion on $\De_{m_0}$ with
dilatation $K=K(\de(a,\partial \De_{m_0}))$, which depends on the
conformal distance $\de(a,\partial \De_{m_0})$ from $a$ to the
boundary of $\De_{m_0}$. Let us call this extended motion
$\ti{h}_{m_0}$.

Now, lift the motion $\ti{h}_{m_0}$ via the unbranched covering maps
$R_{a}^{m_j-m_0}$ for $a \in \De_{m_j}$. We get a holomorphic motion
\[
\ti{h}_{m_j}: \De_{m_j} \times A_{m_j}^{a_0} \longrightarrow \hat{\C},
\]
where $A_{m_j}^{a_0} = P_{m_j}(-a_0) \sm P_{m_j+1}(-a_0)$ is an
annulus surrounding the critical value (the $A_{m_j}$ are children
to $A_{m_0}$). Since holomorphic composition does not change the
dilatation, it follows that this lifted motion has the same
dilatation $K$ as $\ti{h}_{m_0}$. Moreover, the annuli
$A_{m_j}^{a_0}$ move holomorphically; set
$A_{m_j}^a=\ti{h}_{m_j}(A_{m_j}^{a_0})$. In other words, $A_{m_j}^a
= D_{m_j}^a \sm D_{m_j+1}^a$.

By Lemma \ref{diff}, we have that $a \in \De_{m_j} \sm \De_{m_j+1}$ if
and only if $-a \in A_{m_j}^a$.

Define the {\em parameter annuli} $\AA_n = \De_n \sm \De_{n+1}$.

\comm{

Take any point $a \in \partial \AA_N$. Thus $a$ belongs to the closure of a
parabubble $P$. If it is in the interior of $P$ then there is a
smallest $k$ such that
$R^k(-a) \in A_{\infty}$ and there is a unique B\"ottcher coordinate
$\psi_a(R_a^k(-a))$ in $\hat{\C} \sm \D$ for $a$ (see Section
\ref{parabubble}). Now let
\[
F: a \mapsto R_{a_0}^{-k} \circ \psi_{a_0}^{-1} \circ \psi_a \circ
R_a^k(-a).
\]
By the uniqueness of the B\"ottcher coordinates, it is clear that $F$
is a homeomorphism from $\partial \AA_N$ onto $\partial A_N^{a_0}$.

Recall that the holomorphic motion $h=h_N: \De \times \hat{\C} \raw
\hat{\C}$. For fixed $z \in B$, for some bubble $B$ in the boundary of
the puzzle of depth $d(N+1)$, and some $a \in
\De$, let $k$ be
smallest such that $R_a^k(z) \in A_{\infty}$. Assume also that $B$
lies in the bubble puzzle of depth $d(N+1)$. Then the function
$\ti{h}_a$ coincides with
\[
f_a(z)=R_{a}^{-k} \circ \psi_{a}^{-1} \circ \psi_{a_0} \circ R_{a_0}^k(z),
\]
which is the expression for the holomorphic motion of the boundary of
$A_N^{a_0}$.

}

Fix the number $N=m_j$ from now on and let $\De_{N}=\De$.
Define a map defined on $\De$, by
\[
H=H_N: a \mapsto \ti{h}_a^{-1}(-a).
\]
We see that $H_N: \AA_N \raw A_N^{a_0}$. On the boundary of $\De$ it
is injective, which follows directly from Lemma \ref{boundary}.

The next issue is to show that the map $H_N$ is quasiconformal with a
definite bound on the dilatation independent of $N$. Here the proof
is again the same as in \cite{PR}; let us differentiate the relation
$\ti{h}_{N}^a(H_N(a)) = -a$. Then we get
\[
\oli{\partial} h_{N}^a (H_N(a)) \oli{\partial H_N(a)} + \partial h_{N}^a
\oli{\partial} H_N(a) = 0.
\]
This implies that the Beltrami coefficient $\mu(a) = \oli{\partial} H_N /
\partial H_N$ satisfies
\[
|\mu(a)| = \frac{|\oli{\partial} h_{N}^a(H_N(a))|}{|\partial
 h_{N}^a(H_N(a))|} =
 \frac{K_N-1}{K_N+1} < 1,
\]
where $K_N$ is the dilatation of $h_{N}^a$. However,
if we consider the conformal representation $\chi: \De_N \mapsto \D$,
The $\la$-Lemma implies that
\[
K_N = \frac{1+|\chi(a)|}{1-|\chi(a)|}.
\]
Since the sets $\De_{m_j}$ is compactly contained in $\De_{m_0}$ for
$j \geq 2$, we get that $K_{m_j} \leq K$, for all $ j \geq 2$.

We claim that the map $H_N$ is injective.
First of all, it is injective on the boundary of $\AA_N$.
Moreover, if we solve the Beltrami equation for $\mu$,
then we get a quasiconformal map $\phi:\AA_N \raw \phi(\AA_N)$, so
that $\oli{\partial} \phi = \mu \partial \phi$. It follows that $H_N
\circ \phi^{-1}$ is conformal. By the Riemann- Hurwitz formula, there
can not be any branch points in $\AA_N$. Since $H_N$ is injective on
the boundary of $\AA_N$, it follows that $H_N \circ \phi^{-1}$
maps $\phi(\AA_N)$ conformally onto $A_N^{a_0}$. It follows that $H_N$
must be a homeomorphism.

Since the annulus $A_{m_0}(-1)$ may be degenerate, we again consider
the complementary annuli $\alpha_m (-1)$.


It follows that
\[
\frac{1}{K} \mod \alpha_{m_j}^{a_0} \leq \mod \tl\alpha_{m_j} \leq \frac{1}{K}
\mod \alpha_{m_j}^{a_0},
\]
where $\tl\alpha_m$ denotes a complementary annulus in the parameter plane.
Since
$$\sum \mod \alpha_N = \infty\text{ we have }\sum \mod \tl \alpha_N = \infty,$$
and we conclude from \lemref{moduli diverge} that the parameter pieces $\De_N$ shrink to a single
point, which has to be $a_0$.

In the non-recurrent case, consider the puzzle of depth $N$ so that
the critical puzzle piece $P_N(-1)$ is disjoint from the postcrtical
set. As the critical value $-a$ varies through $\De_N$ the puzzle at
depth $N+1$ moves holomorphically as in Lemma \ref{initial}. Hence
every annulus $A_N(z)$ moves holomorphically. Extend this
holomorphic motion by Slodkowski's Theorem and denote the extended
motion by $\ti{h}$ similar to the above argument. Since every
annulus $A_n(-a_0)$, for $n > N$, is a univalent pullback of some
$A_N(z)$ (since $R-{a_0}$ is non-recurrent) we can lift the
holomorphic motion $\ti{h}$ to the parameter piece $\De_n$ over
$P_n(-a_0)$. Define a map $H_n: \AA_n \mapsto A_n(-a_0)$ in exactly
the same way as above. The proof of the fact that the parameter
annuli shrink to a single point is now similar to the recurrent case
and we leave the details to the reader.

\comm{ In the non-recurrent case, we can transfer to the parameter
plane a quasiconformal copy of the annulus $P^{a_0}_n \setminus
P^{a_0}_{n+m}$ for an arbitrarily large $m$ and fixed $n$. This
immediately implies the claim. We leave the details to the reader in
this case.}

\end{proof}

\noindent
We conclude:

\medskip
\noindent {\bf Main Theorem, the uniqueness part.} {\it The mating
in Main Theorem is unique.}

\bibliographystyle{plain}


\end{document}